\documentclass[final]{siamltex}
\usepackage{graphicx}
\usepackage{amsfonts, amsmath}
\usepackage{hyperref}
\usepackage{subfig}
\usepackage{color}
\usepackage[rgb]{xcolor}

\hypersetup{pdfborder={0 0 0}}

\begin{document}

\title{Dynamics of the Budyko Energy Balance Model }

\thanks{The author was supported in part by a grant to the University of Arizona from HHMI (52006942). \\ The research described here has been supported by NSF grant DMS-0940363.} 

\author{Esther R. Widiasih}

\maketitle
\begin{abstract}
A dynamical system derived from a conceptual climate model is investigated. The model is often referred to as Budyko-Sellers type Energy Balance Model, or in short, the Budyko model. The focus of this study is on the mathematical formulation that further develops the model by capturing the ice line dynamics, and thereby places the Budyko model in dynamical systems framework. The classical Budyko EBM  is represented  here as an integro-difference equation, then  coupled with an ice line equation. The resulting infinite dimensional system  is shown to possess an attracting one dimensional invariant manifold. The novel formulation captures solutions that others have previously obtained  while also making the formulation and the analysis of the model more precise.   
\end{abstract}

\section{Introduction}
We investigate a dynamical system derived from a conceptual climate model. The model is a version of what is known in climate science as the Budyko Sellers type  Energy Balance Model (EBM). The result of this study is purely mathematical with the main focus on the mathematical formulation that further develops the model and places it in the framework of dynamical systems. In particular, we represent the classical Budyko model as an integro-difference equation, then couple it with an auxilliary equation that prescribes the ice line dynamics. We analyze the stability of the coupled system and we prove that, under certain conditions, the infinite dimensional system has a one dimensional attracting invariant manifold.\\

Budyko-Sellers type EBMs model the annual average temperature distribution, and rest on the idea that the energy received is balanced by the energy emitted by the planet. Any imbalance in the system will result in a temperature change. The state variable of a Budyko EBM is a latitudinally averaged temperature distribution, often referred to as a \emph{temperature profile}. One feedback admitted by a Budyko-Sellers type model is the ice albedo feedback, which is a positive feedback caused by the contrast between the reflectivity of ice cover and ocean.  Much effort has been devoted to the steady state analysis of differential equations based on the two conceptual models pioneered independently by Mikhail Budyko \cite{bud69} and William Sellers \cite{sellers69} in late sixties eg. \cite{n79, n75, smg88, heldsuarez}. Here, we are concerned with the formulation of the model in order to clarify the role of the ice albedo feedback and make precise the dynamics of both the transient and the steady states. \\

Geologic evidence shows that the northern ice cap had in the past advanced to near tropical regions \cite{mac}. In the more recent past, during the last ice age around 12,000 years ago, 2 mile thick glaciers advanced as far as Illinois  and Wisconsin \cite{ehlers}. Our planet has also experienced warmer climates with both poles having small to no ice cap. The movement of the ice extent or ice cap boundary as a transient response of the climate motivates our study in this work. Indeed, the solar variations of the past few million years are too small to account for the large fluctuations in climate such as glaciation/ interglaciation events. Therefore, climate feedbacks, like those due to greenhouse gases and ice albedo feedback, provide explanations for these major changes and understanding the transient climate response to the feedbacks is of central concern in climate science. \\

At the heart of the Budyko model is the ice albedo feedback. The reflectivity contrast between a dark blue ocean and a shiny  white ice covered surface causes a positive feedback in a climate system. One may imagine an aquaplanet consisting only of water and ice. A warming climate in this planet causes ice to melt, hence the retreat of an ice cap, which causes an increase exposure of the ocean surface, more energy absorption, and thus a warmer climate. A similar positive feedback also works in the opposite direction: a cooler climate causes a favorable condition for ice to grow, hence the advance of an ice cap, which causes more energy reflected to space, and thus a cooler climate. In such a system, one wonders if extreme states (ie. a completely ice free or a completely ice covered state) are the only equilibria of the system. A Budyko-Sellers type model is a simple climate model that attempts to answer this question based on  the energy balance principle. Here, we focus on an EBM similar to the one introduced by Budyko. \\

The goal of this paper is to further develop the Budyko EBM with the ice albedo feedback to place it within a dynamical system framework. Therefore, we focus on the treatment of the transient dynamics of the ice cap boundary location or what we call here the \emph{ice line}. We obtain an equation dictating the ice line movement by exploiting the ice line condition prevalent in current literature. This condition is: when the system is at equilibrium, the  annual average temperature at the ice line location must be the prescribed critical temperature. The conventional approach, which only considers Budyko's EBM for the temperature profile without a coupling to an ice line dynamics, first finds the equilibrium of the temperature profile, and thereby implicitly assumes that such equilibrium exists and is attracting in some function space. It then selects an ice line  location using some prescribed conditions, and finally, performs the stability analysis on the selected ice lines relative to some parameter changes.  In contrast to the established approach, in our formulation we first introduce an ice line equation to be coupled to the Budyko EBM equation of the temperature profile, then show the existence of an attracting one dimensional invariant manifold for the coupled system. The ice line dynamics, ie. the advance and retreat of the ice cap, then takes place on this one dimensional invariant manifold, and in this way, the novel formulation reframes the ice line selection mechanism in the context of stability.  The bifurcation analysis as parameters vary is treated once the stability of the equilibrium states is established.\\   


The paper is organized as follows. Section \ref{back} treats the background of the Budyko Energy Balance Model, and discusses the equilibrium solutions and introduces the novel  formulation as the Budyko equation is coupled with an ice line equation to form a fast-slow system.  Section \ref{pf} states and proves the results: the existence of a one dimensional attracting invariant manifold and the stability analysis.  Section \ref{anim} presents some bifurcation diagrams and animations that illustrate the transient dynamics of the fast slow system and some discussions about the simulation.  Section \ref{dir} concludes the article and discusses some future direction. \\

\section{Background}\label{back}

\subsection{Budyko's Original Concept}

The energy balance at a given latitude on the planet depends on three main aspects of the climate: the incoming solar radiation (insolation), the planet's outgoing radiation and the heat transported to and from that latitude. Budyko's original paper \cite{bud69} has the following two main relationships on these three aspects of the climate:\\

\begin{enumerate}
\item  The equation of the heat balance of the Earth-atmosphere system for mean annual condition
\begin{align}
\text{gain or loss of heat} &= Q (1-\alpha)- I\\
&=Q(1-\alpha)-(A+BT)
\end{align}

where the gain or loss of heat comes \lq\lq as a result of atmosphere and hydrosphere circulation, including heat redistribution of phase water formations \rq\rq (from p. 613 of \cite{bud69}, note that our notations are slightly different from Budyko's original article). Here, 

\begin{tabular}{l p{10cm}}
\hspace{0.15cm}& $I = A+BT$ is the outgoing radiation (at a given latitude)\\ 
 \hspace{0.15cm}& $T$ is the annual mean temperature at a given latitude,  called the \textit{temperature profile} \\
 \hspace{0.15cm}&$Q$ is the solar radiation  coming to the outer boundary of the atmosphere\\
 \hspace{0.15cm}& $\alpha$ is the albedo of the surface.\\
\end{tabular} \\

\item The horizontal heat transfer distribution 
\begin{equation}
\text{gain or loss of heat} = C ( T-\overline{T} ).
\end{equation}

 \begin{tabular}{l p{10cm}}
 Here, & $\overline{T}$ is the globally-averaged annual mean temperature\\
 \hspace{0.15cm}& $C$ is a proportionality constant.\\
\end{tabular}\\

\item When the system is in balance, there is no gain or loss of energy (heat)  at each latitude and the following equation holds: 

\begin{equation} \label{mainbud}
0=Q(1-\alpha)-I-C(T-\overline{T}).
\end{equation}

\end{enumerate}

A similar concept was proposed using a different formulation by William Seller around the same time \cite{sellers69}. Based on this fundamental relationship and by treating the temperature profile as a static system, Budyko calculated and analyzed the effect of exogenous changes in solar forcing (due to variations in Earth's orbital elements) and atmospheric transparency (due to volcanic dust). In a more recent work,  Tung \cite{tung07} established the modern parametrizations of the Budyko model and included a nice exposition of the model. In the most recent work McGehee and Lehman \cite{mcgleh} combined a Budyko's model version formulated in \cite{tung07} with Earth's orbital elements and compared it with the climate data. \\

The state variable of the Budyko model \eqref{mainbud} is the temperature profile. Each temperature profile is therefore a function in some infinite dimensional space. To determine  the ice line,  the conventional approach first finds the steady state temperature profile, then uses a \emph{selection mechanism} specifying a critical temperature to determine the ice line location (see eg. \cite{tung07, hsu}).  The latter article proceeded further by introducing an ice line condition based on the ice-water boundary temperature, without explicitly stating the equation that governs the dynamics of the location of this ice-water boundary. \\

Many authors have treated the fundamental relationship \eqref{mainbud} as the equilibrium of some dynamical system while varying the parameters of the model. Using this approach, the works by North, et al in \cite{n79, n75, smg88, heldsuarez}  studied the stability of a similar system with a diffusive transport, while the work by  Hsieh and Su in \cite{hsu} studied the stability by including a time variation of the ice boundary. The stability analyses done by these authors were based on the perturbation of some parameters near the equilibrium temperature profile, and the reader will find a more detailed discussion about this in Section \ref{eqtempf}.  However, most authors have not suggested a particular function space for the temperature profile. Furthermore, since most authors placed more emphasis on the application and the parametrization of the model, none suggested an explicit description of the ice line transient dynamic as a response to temperature changes. In order to have a complete and a clear picture of the ice albedo feedback concept and to analyze the stability of the system, we argue that it is necessary to address both the issues of the function space and the ice line dynamics. Once complimented with an appropriate function space and the ice line dynamics, the Budyko model can be treated as a mathematical object, and in particular, as a dynamical system for which a rigorous stability analysis can be developed. \\

In the following subsections, we will lay out the description of each aspect of this dynamical system. The discussion will address:
\begin{enumerate}

\item Budyko's EBM as a difference equation
	\begin{enumerate}
	\item The state variable $T$
	\item The solar forcing distribution
	\item The outgoing solar radiation
	\item The transport term
	\item The temperature profile equilibrium solutions
	\end{enumerate}

\item The Ice line equation\\

\end{enumerate}

\subsection{Budyko's EBM as a Difference Equation}

Following Budyko \cite{bud69}, the change in the annual average temperature in zone $y$ from year $n$ to year $n+1$ is:

\begin{align} \label{mapbud}
\notag  \Delta T&=f(T,\eta):= \frac{1}{R} \left( Qs(y)(1-\alpha(\eta, y))-(A+B T(n,y))-C\left (T(n,y)-\int_0^1 T(n,y)dy \right) \right)\\
\end{align}

\noindent with $\Delta T = T(n+1, y)-T(n,y)$. \\

In what follows, we describe each term of equation \eqref{mapbud}.\\

\subsubsection{The  state variable $T$} 
As explained by North in \cite{n75}, the latitude $\theta$ is more effectively represented by the location $y$, where $$y=\sin(\theta).$$ In this representation, the infinitesimal band $$dy=\cos(\theta)d\theta$$ is proportional to the area of the latitudinal strip parallel to $y$.  Here, we assume the symmetry of the planet, and so $0 \le \theta \le 90^o$ or $0 \le y \le 1$.\\

The state variable $T$ in Budyko's EBM is the average temperature profile during year $n$.  The temperature profile $T(y)$ is measured in Celsius and is a function
of the independent variable $y$.  Thus, $$T = T(n,y)$$ is the surface temperature during year $n$ at latitude $\arcsin(y)$ averaged over the entire year and averaged over all longitudes.  The constant $R$ is the heat capacity of the surface of the planet.\\

We remark here that in his 1969 paper \cite{bud69}, Budyko formulated the change in the temperature profile $T(y)$ in yearly increments. All other quantities, eg. $Q$, $s(y)$, $\alpha(\eta,y)$, etc. in the  equation \eqref{mapbud}  are also yearly averages, having small variations from year to year.  Other authors, eg.  \cite{n75, n84, n90, hsu, heldsuarez} have treated a similar model in the continuous time framework, where the state variable is $T(t,y)$ which depends on the continuous time variable $t\in \mathbb{R}$. In the continuous time framework, the Budyko equation \eqref{mapbud} may be written as a differential equation
$$\frac{\partial T}{\partial t}= \frac{1}{R} \left( Qs(y)(1-\alpha(\eta, y))-(A+B T(t,y))-C\left (T(t,y)-\int_0^1 T(t,y)dy \right) \right)$$
\noindent The reader may observe that both the discrete and the continuous formulations of the Budyko model are both approximations of Earth's climate. Furthermore, in both  frameworks, the systems will retain the same key features eg. the same equilibrium solution, bifurcation diagrams, etc. Therefore, it stands to reason that similar arguments for the stability analysis  apply and that the main result of this article would also hold in the continuous time framework.

\subsubsection{Incoming energy $Qs(y)(1-\alpha(\eta,y))$}
The constant $Q$ is the solar constant, representing the amount of energy received from the sun at the top of the atmosphere. The function $s(y)$ is the latitudinal distribution of that energy which can be computed from earth's orbital elements. Many authors,  such as Tung \cite{tung07}, and North \cite{n75} used the Legendre polynomial function approximation $$s(y)=1-0.482\frac{3y^2-1}{2}$$ for this distribution. We will use the same approximation as well. \\

The fraction of the radiative energy absorbed by the planet at $y$ given that the ice line is at $\eta$ is represented by the factor $$1-\alpha(\eta,y).$$ We will often write $\alpha(\eta,y)$ as $\alpha(\eta)(y)$ to emphasize the fact that for each $\eta$, $\alpha(\eta)(y)$ is a real valued function over the $y-$interval.\\

Here, it is assumed that the surface is either water (ice free) or ice covered and that there is only one ice line, $\eta$.  Since ice reflects more sunlight than  ocean does, the surface covered with ice has higher albedo than that covered with ocean. The albedo function $\alpha(\eta,y)$ we use in this paper is smooth and is ice line dependent.\\

\begin{definition} {\textbf{ Smooth, Ice line-Dependent Albedo}}\\
Given that the ice line is at $\eta$, the albedo at $y$ is
$$\alpha(\eta,y)=\frac{\alpha_1+\alpha_2}{2}+ \frac{\alpha_2-\alpha_1}{2} \cdot \tanh[M \cdot (y -\eta)]$$
where $\alpha_1= \text{ the ocean or water albedo }$ and $ \quad \alpha_2 = \text{ the ice albedo.} $\\
\end{definition}

The parameter $M$ represents the steepness of the albedo function near the ice line and is a fixed quantity. The albedo at the ice line is always the average of the ice and water or ocean albedos. Equatorward of the ice line, the albedo is approximately  the ocean albedo while poleward it is the ice albedo. Notice that both $\alpha_1$ and $\alpha_2$ are fractions of energy reflected, hence dimensionless. Since the reflectivity of ice is higher, the value of $\alpha_2$ is larger than $\alpha_1$.  A similar smooth albedo formulation is used in \cite{dg81}, while a piecewise function is used by others, eg. Tung \cite{tung07}. In the spirit of the piecewise albedo function used in the latter work, we chose $M \gg 1$. \\

\subsubsection{Outgoing longwave radiation}

The absorbed radiative energy contained in the first term of the left hand side of equation \eqref{mapbud} is balanced by the re-emission or re-radiation term and the transport term. The re-emission term at $y$ is  $$-\left( A+BT(y) \right)$$  which is an  approximation of both the Stephan-Boltzman's law of black body radiation and the greenhouse gas effect on the atmosphere. An alternative parametrization, based on the work by Graves and North \cite{graves93}, was used by Tung in the Budyko model formulation in \cite{tung07}, where the constants $A$ and $B$ were derived from fitting a linear function through satelite data of Outgoing Longwave Radiation (OLR) at the top of the atmosphere.  \\

\subsubsection{Transport}\label{transportsubsec}

The energy transport is represented by the term, 
\begin{equation}\label{transport} C \cdot \left( T(y)-\int_0^1 T(y)dy \right).\end{equation}
The constant $C$ in the transport term is traditionally chosen so that one of the equilibria fits the current climate with an ice line near the pole. The transport term presented here, which Engler and Kaper refered to as the heat \emph{transfer} \cite{kaper}, is an alternative the diffusion transport. In \cite{n90}, North studied the Budyko model with diffusion transport term, and  specified a boundary condition that no heat flux shall pass through the equator and the pole. Observe that the diffusive transport is inherently a local process which explicitly specifies the heat transport process in a neighborhood of each zone $y$, hence the boundary condition at the equator and the pole.  In contrast, the transfer transport  \eqref{transport} is a global process that does not necessitate a special treatment of the heat flux at the boundaries, ie. the equator, for $y=0$ or pole $y=1$. 

\subsection{The temperature profile equilibrium solutions} \label{eqtempf}
The equilibrium solution of equation \eqref{mapbud}  can be solved by fixing $\eta$. For a fixed $\eta$, one first integrates both sides of the equation $f(T,\eta)=0$ in equation \eqref{mapbud} with respect to $y$ over the unit interval, then solves for the functional $\overline{T}^*(\eta):=\int_0^1 T^*(\eta)(y)dy$. After reinserting the $\eta$-dependent expression for $\overline{T}^*(\eta)$ in the equation $f(T,\eta)=0$, one then solves for $T^*(\eta)$. We obtain 
\begin{equation} \label{meaneqtemp01}\overline{T}^*(\eta)=\frac{Q\overline{\alpha}(\eta)-A}{B} \end{equation}

where $\overline{\alpha}(\eta)=\int_0^1 s(y)(1-\alpha(\eta)(y))dy$. Explicitly, the equilibrium temperature profile of equation \eqref{mapbud}, with parameters as in table \ref{parametertable} is:
\begin{equation}\label{eqtemp01}
T^*(\eta)(y)=\frac{Q \cdot s(y) \cdot(1-\alpha(\eta)(y))-A+ \frac{C}{B}\left(Q\overline{\alpha}(\eta)-A \right)}{B+C}
\end{equation}

\noindent The figures on the bottom row of Figure (\ref{icf}) are three examples of local temperature profiles for $\eta = 0.1, 0.5, 1$, when $y$ is restricted to the unit interval. \\

A \emph{selection mechanism} of the ice line prevalent in the literature is the the following: for an ice line location $\eta$ to be an ice line equilibrium, the  temperature at that location must reach the critical temperature, $T_c$. Therefore, the ice line equilibria is the solution of the equation 
\begin{equation} \label{icelineselection}
T^*(\eta)(\eta)=T_c.
\end{equation} 
 Using the parameters in \cite{tung07}, one easily checks that there are two ice lines that satisfies equation  \eqref{icelineselection}, $\eta_1 \cong 0.2$,  and, $\eta_2 \cong 0.95$. It is customary to call the equilibrium $\eta_1$ \textbf{the big ice cap} and $\eta_2$ \textbf{the small ice cap}. The first equilibrium will be shown to be a saddle, while the second is a sink. In additional, the two physically attainable state that we must consider are \textbf{snowball} $\eta=0$ and \textbf{ice free} $\eta=1$. \\

One may notice that a curve of a parameter, say $Q$, as a function of the ice line  $\eta$, can be obtained from the ice line selection mechanism, equation \eqref{icelineselection}. Using intuitive arguments, Budyko in 1972 argued that the branch where $\frac{dQ}{d\eta}>0$ is stable while the branch where $\frac{dQ}{d\eta}>0$ is unstable  \cite{bud72, heldsuarez, tung07}. Held and Suarez  extended this intuitive argument with a stability analysis based on a linear perturbation on the $\overline{T}^*$. (See details on p. 615-616 in \cite{heldsuarez}). Notice that these stability analyses operates in a one dimensional space despite the fact that the temperature profile in the Budyko model is a function over the unit interval, hence, an infinite dimensional object. Such result points to the existence of a one dimensional attracting invariant manifold. In what follows, we expand on the  work previously done by specifiying a function space for the temperature profile and introduce an ice line equation, thereby providing a context for which  a one dimensional attracting invariant manifold exists under the conditions set by the Budyko model and  the ice line selection mechanism. \\

\subsection{The Ice line Equation}
The equation for the ice line movement must reflect the following: \\

\begin{tabular}{l p{10cm}}
 \hspace{0.18cm} & 1. Poleward or equatorward ice line movement when the ice line is, respectively, too hot or too cold with respect to a reference critical temperature. \\
 \hspace{0.18cm}&  2. The ice line moves at a much slower speed than the atmospheric temperature reflecting the rate at which glaciers advance or retreat.\\
 \hspace{0.18cm}& 3. When the ice line equation is coupled with equation \eqref{mapbud}, the resulting equilibria must be consistent with the conventional approach.\\
\end{tabular}\\

\noindent We propose the following equation to account for the one year zonal change of the ice line $\eta$:
\begin{equation}\label{iceeq}
\Delta \eta = g(T,\eta):=\epsilon(T(\eta)-T_c)
\end{equation}

\noindent with $\Delta \eta = \eta(n+1)-\eta(n)$,  $T_c$ a reference critical temperature,  and $\epsilon \ll 1$. \\

That is, the ice line moves at a (much slower) rate proportional to the difference between the ice line temperature and the reference critical temperature. Some authors used $T_c=-10$ for a modern day parametrization, eg. \cite{tung07}. \\

Assuming a certain thickness of ice,  the value of $\epsilon$ can be computed to be in the order of $10^{-5} \left( /^oC\text{ year} \right)$ (see \cite{mcgw}). Here, the value  of $\epsilon=0.01$ is used to create the dynamic version of the simulation  in Figure \ref{dyn}. This value is chosen so that the animation run time is reasonably short.\\

Although all of  the aforementioned authors in the previous studies indicate a critical temperature at the ice line, none had mentioned the ice line equation explicitly as written in equation \eqref{iceeq}. Some authors indicate an ice line rate of change as  byproducts of the stability analysis as indicated in the discussion in subsection \ref{eqtempf}, eg $\frac{d\theta_0}{dt}$ in the Appendix A of \cite{heldsuarez} and $\frac{dx_s}{dt}$ on page 2274 of \cite{hsu}, however, an explicit formulation of the ice line transient dynamics  has been absent.\\

\section{Statements and Proofs of the Main Result and its Corollaries}\label{pf}
\subsection{The Main Result} The object of our analysis is the following fast-slow system of integro-difference equations consists of the equations from \eqref{mapbud} and \eqref{iceeq}:
\begin{align}\label{budt} 
\Delta T &= f(T,\eta) :=  \frac{1}{R} \left[ Qs(y)(1-\alpha(\eta,y))-(A+B T)-C\left(T-\int_0^1 Tdy \right) \right]  \\ 
\Delta \eta &= g(T, \eta) := \epsilon(T(\eta)-T_c) \notag
\end{align}

where $\eta, y \in [0,1].$\\

The following definition defines the map based on system \eqref{budt}:\\

\begin{definition} {The Budyko Map} \\
The Budyko map $\widehat{m}$ is the map associated with system \eqref{budt}
$$\widehat{m}([T,\eta]) = \left[ T+f([T,\eta]), \eta+g([T,\eta]) \right].$$
\noindent Denote by $\widehat{m}^k$ the $k-th$ iteration of the map $\widehat{m}$.\\
\end{definition}

\noindent Observe that if system (3.2) were formulated instead as a differential equation, the Budyko map would be the forward Euler approximation with time step one year.\\

Some conditions on the Budyko model's parameters are necessary to show the main result.  We collect the conditions in the following definition: \\

\begin{definition} {The Budyko Parameters Condition} \\
\noindent We say that the system \eqref{budt} satisfies the Budyko parameters condition if:
\begin{description}
\item{1.} The parameters $Q,A,B,C,R,M$ are all positive,\\ with $Q \gg B \text{, } A \gg B \text{, } Q \gg C \text{,  } A \gg C \text{ and } M \gg 1$.\\
\item{2.} $0<1-\frac{1}{R}(B+C) < 1$, or similarly, $0<\frac{1}{R}(B+C)< 1$ .\\
\end{description}
\end{definition}

\noindent These conditions are physically reasonable and are consistent with the parameters used by other authors, eg. \cite{tung07}. The second condition, in particular, is necessary to show the existence of an invariant manifold.\\

We now discuss the function space of the temperature profile. It is not unreasonable to demand each temperature profile to be bounded, since physically, planetary temperature is always bounded. Denote by $\mathcal{B}_{I}$ \textbf{the space of the  continuous function} over the unit interval $I=[0,1]$ with the $sup$ norm. The phase space for the temperature profile and ice line we will consider is $\mathcal{B}_I \times I$, with the norm on the product space $\| \cdot \|_{{{\mathcal B}_I} \times I} = max \left\lbrace \| \cdot \|_{{\mathcal B}_I}, |\cdot| \right\rbrace$. The main result of this paper is the following theorem:\\
\begin{theorem}\label{mainresult} \textbf{Main result} Suppose that $\eta, y \in  I=[0,1] $ and the system \eqref{budt} satisfies the Budyko parameter condition. Then for all $\epsilon$ sufficiently small:
\begin{enumerate}
\item \textbf{Attractivity.} Under the Budyko map $\widehat{m}$, there exists a one dimensional locally exponentially attracting Lipschitz continuous manifold with boundary $$\mathcal{M}_I:=\left\lbrace  [\phi^*(\eta)(y), \eta]: \eta, y \in I \right\rbrace \subset  \mathcal{B}_I \times I$$

\item \textbf{Invariant alternatives.} The action of the Budyko map on $\mathcal{M}_{I}$ is either one of the following:\\
\begin{enumerate} 
\item $\widehat{m}([\phi^*(\eta)(y), \eta]) = [\phi^*(\xi)(y), \xi]$ , with $\xi \in [0,1]$\\
\item[] OR\\
\item  $\widehat{m}([\phi^*(\eta)(y), \eta]) = [T(y), \xi]$, $T\in \mathcal{B}_I$ and $\xi \notin [0, 1]$\\
\end{enumerate}
\end{enumerate}
\end{theorem}

\noindent The following two corollaries are pertinent to the physical interpretation of the model. \\

\begin{corollary} \label{critset}
The invariant manifold $\mathcal{M}_I$ is within $O(\epsilon)$ of the critical set $$\mathcal T_I^*:=\left\lbrace [T^*(\eta)(y),\eta]: \eta, y \in I=[0,1] \right\rbrace$$ with $T^*(\eta)(y)$ as defined in \eqref{eqtemp01}. \\ 
\end{corollary}

\begin{corollary} \label{stability}
With the parameters as in the table (\ref{parametertable}), the ice free Earth and the big ice cap are unstable while the small ice cap and the snowball earth are stable.\\
\end{corollary}

\begin{proof}

In order to simplify the analysis of the ice line dynamics at the end points, that is for $\eta=0$ and $\eta=1$, we embed the ice line interval and the domain of the temperature profile  in the real line. We will also embed the right hand side of \eqref{budt} in such a way that preserves the dynamics inside the closed unit interval, then show that a similar result holds in the embedded space. The reader will observe that this embedding is a mathematical convenience to show the existence of an invariant manifold and to allow  the dynamical analysis of the polar and the equatorial ice line, particularly in the details of the proof of Lemma \ref{preimage}.\\

Denote by $\mathcal{B}$ \textbf{the space of the bounded continuous function} over $\mathbb{R}$ with the $sup$ norm $\| \cdot \|_{\infty}$. The  \emph{temperature profile and ice line} phase space is then chosen to be $\mathcal{B} \times \mathbb{R}$, endowed with a similar product norm as the restricted space $\mathcal{B}_I \times I$. That is, $\| \cdot \|_{\mathcal{B}\times\mathbb{R}} := max \left\lbrace \| \cdot \|_{\mathcal{B}}, |\cdot |_{\mathbb{R}} \right\rbrace $. Note that one may continuously embed any function $f$ in $\mathcal{B}_I$ into $\mathcal{B}$. On the other hand, any function $F \in \mathcal{B}$ can be restricted to $\mathcal{B}_I$. Furthermore, regardless of the embedding, the norm of the embedded function majorizes its restricted counterpart. We then  proceed to prove that the attractivity and invariant results holds in the embedded space $\mathcal{B} \times {\mathbb R}$. The main result Theorem \ref{mainresult} immediately follows.\\


The following definitions give the embedding of the system \eqref{budt} into $\mathcal{B}_{\mathbb{R}} \times \mathbb{R}$, the extended Budyko map $m$ of the Budyko map $\widehat{m}$ based on this embedding, and the extended critical set $\mathcal{T^*}$ of the set $\mathcal{T}^*_I$.  

\begin{definition} \label{budxtdef}
 \begin{align} \label{budxt}
&\Delta T  := F \left( \left[ T,\eta \right] \right) \\
&\Delta \eta := G \left( \left[ T,\eta \right] \right) \notag
\end{align}
   
where we define:  \\

\begin{equation} \label{iceeqxt} G([T,\eta]) := \epsilon [T(\eta)-T_c] \end{equation}

and \\

$F([T,\eta])$:=
\begin{equation} \label{budembed}
 \begin{cases}
 & \frac{1}{R} \left( Q\cdot s(0) \cdot [1-\alpha(\eta)(y)] - [A+BT(y)] + C [\overline{T}-T(y)] \right)  \text{,  when } y <0\\
& \frac{1}{R} \left( Q\cdot s(y) \cdot [1-\alpha(\eta)(y)] - [A+BT(y)] +C  [\overline{T}-T(y)  ]\right) \text{,  when } 0 \leq y  \leq 1\\
& \frac{1}{R} \left( Q\cdot s(1) \cdot [1-\alpha(\eta)(y)] - [A+BT (y)] +C [\overline{T} -T(y)]\right) \text{,  when } 1<y.\\
\end{cases}
\end{equation}

with $\overline{T}= \int_0^1 T(y) dy$.\\
\end{definition}

The reader may be tempted to interpret the embedded Budyko equation \eqref{budembed} as a treatment of the temperature profile at the boundaries. ie. the equator and the pole. Again, we emphasize that the embedding the Budyko equation is a mathematical technique to simplify the analysis, and is not a substitute for the physical modeling of the boundary condition for the temperature profile $T(y)$ in  the Budyko model, akin to the diffusive transport version discussed in \S  \ref{transportsubsec}. Indeed, the model may be improved further and the reader will find such discussion in \S \ref{futwk},  Future Directions.\\

\begin{definition} {The Extended Budyko Map} \\
\noindent The Budyko map $m$ is the  map associated with the difference equation \eqref{budxt}:
ie. when the initial conditions are $T_0=T(0,y)$ and $\eta_0=\eta(0)$, 
$$m([T_0,\eta_0]) = \left[ T_0+F([T_0,\eta_0]), \eta_0+G([T_0,\eta_0]) \right].$$
 \\
\noindent Denote by $m^k$ the $k$-th iteration of the map $m$.\\
\end{definition}

\begin{definition}{The Extended Critical Set}\\
The critical set of the extended system \eqref{budxt} is the subset of $\mathcal{B}$ $$ \mathcal{T^*}:= \left\lbrace T^*(\eta)( y): \eta, y \in \mathbb{R}  \right\rbrace $$
where for a fixed $\eta$,  

\begin{equation*}
T^*(\eta)( y)=
\begin{cases}
&\frac{Q \cdot s(0) \cdot(1-\alpha(\eta)(y))-A+C\int_0^1{T^*(\eta)(y)dy}}{B+C} \text { ,  when } y<0\\
&\frac{Q \cdot s(y) \cdot(1-\alpha(\eta)(y))-A+C\int_0^1{T^*(\eta(y))dy}}{B+C} \text{ ,  when } 0 \leq y \leq 1\\
&\frac{Q \cdot s(1) \cdot(1-\alpha(\eta)(y))-A+C\int_0^1{T^*(\eta)(y)dy}}{B+C} \text{ ,  when } 1<y.\\
\end{cases}
\end{equation*}
\end{definition}

\noindent Notice that $\mathcal{T}^* \subset \mathcal{B} \times \mathbb R$ and is the set of temperature profiles $T(y)$ that solve the equation $F([(T(y), \eta])=0$.

One may view  $T^*(\eta)(y)$ as the mapping from the reals to the function space  $\mathcal{B}$ where each ice line $\eta$ is mapped to a critical set $T^*(\eta)(y)$. In this way, $T^*$ acts as a graph mapping $\eta \in \mathbb{R}$ into the set of bounded continuous functions $\mathcal{B}$.  \\

We extend this notion of graph over the real line in a more general sense:\\

\begin{definition} \label{graphspace}  \textbf{The space of graphs} \\
We denote by $\mathcal{G}$ the set of graphs:
$$\mathcal{G}:= \left\lbrace \Phi:\mathbb{R} \rightarrow \mathcal{B} | \Phi \text{ is uniformly bounded and continuous over } \mathbb{R}  \right\rbrace $$
endowed with the uniform norm $$\|\cdot \|_{\mathcal{G}} = sup_{\eta \in \mathbb{R}} \| \cdot\|_{\mathcal{B}}.$$
\end{definition}

Note here that since $\mathcal{B}$ is a Banach space, hence complete, and since $\mathbb{R}$ is a metric space, then $\left\lbrace \mathcal{G}, \| \cdot \|_{\mathcal{G}} \right\rbrace$ is a complete metric space. Also, observe that the albedo function $\alpha$ is a graph in $\mathcal{G}$, where for each ice line $\eta$, $\alpha(\eta)( y)$ is smooth, hence Lipschitz, with a Lipchitz constant $M$ and a (sup) norm $\alpha_2$, the ice albedo. This observation will be important in the following characterization of the graph of the equilibrium temperature profiles. \\

The extended critical set will play an important role in the construction of the invariant manifold. Observe that  for each $\eta$, $T^*(\eta) \in \mathcal{B}$ is a Lipschitz function with the steepest decrease at the ice line location $y=\eta$, due to the albedo function $\alpha(\eta)(y)$. It is straightforward to check that for all $\eta$, the Lipschitz constant of the bounded continuous function $T^*(\eta)$ is uniformly bounded.  The ice line temperature function $h(\eta):=T^*(\eta)(\eta)$, viewed as a function of $\eta$, is also a bounded continuous real valued function over the real line, and is also observed to  be Lipschitz continuous with uniformly bounded Lipchitz constant. Based on these observations, we define the constant $L$ that will serve as the upper bound for many of the estimates to come. Our definition of the upper bound constant $L$ here is admittedly coarse, as one may find a smaller upper bound. Nonethelesss, we choose this definition because it simplifies the analysis considerably.\\

\begin{definition}\label{L} \textbf{The constant $L$} \\ We define $$L:=  \max \left\lbrace \|T^* \|_{\mathcal{G}},  Lip_{\mathcal{G}}(T^*), Lip_{\mathcal{B}}(T^*(\eta)), Lip(T^*(\eta)(\eta)).\right\rbrace$$ \\
\end{definition}

where $$Lip_{\mathcal{G}}(\Phi):= \sup_{\eta \ne \zeta; \eta, \zeta \in \mathbb{R}}\frac{\|\Phi(\eta)-\Phi(\zeta) \|_\mathcal{B}}{|\eta - \zeta |}$$
and $$Lip_{\mathcal{B}}(\Phi(\eta)) := \sup_{x \ne y; x \, y \in \mathbb{R} } \frac{|\Phi(\eta)(x)-\Phi(\eta)(y) |}{|x - y |}.$$ 

Although we omit such details here, observe that the constant $L$ may be expressed in terms of the parameters of the Budyko equation. For example, one may show that 
$$|\overline{T}^*(\eta)| \le \frac{(QM+A)(B+C)}{B} $$ and so, $$\|T^*(\eta)\|_{\mathcal{B}} \le QM+A+C(QM+A).$$ Since  $\|T^*(\eta)\|_{\mathcal{B}}$ is independent of $\eta$, then 
$$\|T^*\|_{\mathcal{G}} \le QM+A+C(QM+A).$$

The following result shows the existence of the attracting invariant manifold for the extended system \eqref{budxt}.\\

\begin{theorem} \label{invariant}\\
Suppose that the function $F$ in \eqref{budxt} satisfies the Budyko parameter condition. Then for all $\epsilon$ sufficiently small, there exists an attracting one dimensional invariant manifold for the  map $m$ associated with the Budyko equation \eqref{budxt}. That is,
\begin{description}
\item{1.} There exists a Lipschitz continuous map 
$$\Phi^*:\mathbb{R} \rightarrow \mathcal{B}$$ 

such that $\left[\Phi^*(\mathbb{R}), \mathbb{R} \right]$ is invariant under the Budyko  map $m$, ie. $$m ( [\Phi^*(\mathbb{R}), \mathbb{R}] ) = [\Phi^*(\mathbb{R}), (\mathbb{R})] $$ 

\item{2.} There exists a closed set $\mathcal{C} \subset\mathcal{B}$ in which the invariant manifold is attracting. That is: if $T_0 \in \mathcal{C}$, the distance between the $k$-th iterate of $T_0$ and $\Phi^*(\mathbb{R})$, ie. $$\sup_{\eta \in \mathbb{R}} \| m^{k}([T_0, \eta])-(\Phi^*(\eta), \eta)\|_{\mathcal{B} \times \mathbb{R}}$$ decreases exponentially. \\
\end{description}
\end{theorem}

The proof is based on the Hadamard's Graph Transform method \cite{robinson}. The idea of the proof is to show the existence of a Lipschitz graph over the reals into the Banach space $\mathcal{B}$ that is invariant under the action of the extended Budyko map. This result is possible since for a small enough $\epsilon$, the extended Budyko  map is a contraction on the closed set  $\mathcal{C}$. \\

The statement of  Theorem \ref{invariant} gives us some hints how to proceed with the proof.  We  start by defining the Budyko graph transform operator acting on the space of graphs from  Definition \ref{graphspace}.  The last part of the proof is showing the contraction property of the Budyko graph transform operator.\\

First, consider the following definitions and lemmas:\\

\begin{definition} \textbf{The subset $\mathcal{G}_L$}\\
The set $\mathcal{G}_L$ is defined to be:
$$\mathcal{G}_L:=\left\lbrace \Phi \in \mathcal{G}: \| \Phi \|_{\mathcal{G}} \le L \text{, } Lip_{\mathcal{G}}(\Phi) \le L \text{, and for each } \eta \text{, } Lip_{\mathcal{B}}(\Phi(\eta)) \le L \right\rbrace. $$ 
\end{definition}
 
The following observations are immediate consequences of the definition of $\mathcal{G}_L$ and are necessary to show some upcoming contraction properties: \\

\begin{description} 
\item{1.} The set of graphs $\mathcal{G}_L$ is sequentially  closed in $\mathcal{G}$.\\ 
\item{2.} If $\Phi \in \mathcal{G}_L$ then  for any  real values $x,y$, $$\frac{|\Phi(x)(x)-\Phi(y)(y)|}{|x-y|} \le 2L,$$ and therefore, $Lip(\Phi(x)(x)) \le 2L$. \\
\end{description}

We now show that each ice line has an ice line preimage under $G$. Notice the convenience of the embedding of $[0,1]$ in the real line comes into play here, as one need not fuss with the question whether the preimage of the ice line indeed lies in the domain of the temperature profile.\\
  
\begin{lemma} \textbf{Ice line Preimage } \label{preimage}\\
If $\epsilon<\frac{1}{2L}$, then for each $\eta \in \mathbb{R}$ and each $\Phi \in \mathcal{G}_L$, there exists $\xi \in \mathbb{R}$, with $\xi$  depending on $\eta$ and $\Phi$, such that  
\begin{equation} \label{xi} \eta = G(\Phi(\xi), \xi)=\xi + \epsilon \left( \Phi(\xi)(\xi) - T_c \right). \end{equation} 
\end{lemma}

\begin{proof}
Let $\Phi \in \mathcal{G}_L$ and $\eta \in \mathbb{R}$ be given. We will show that there is a unique function  $k=k_{\Phi} \in \mathcal{B}$, ($k$ depends on $\Phi$), such that for the given $\eta$, the following equality holds:
$$ \eta = \left( \eta + k \left( \eta \right)\right) + \epsilon \left( \Phi \left( \eta + k \left( \eta \right)\right) \left( \eta + k \left( \eta \right)\right) - T_c \right).$$  

Define the transformation $\mathbb{T}$ acting on $k \in \mathcal{B}$ as:
$$ (\mathbb{T}k)(x) = \epsilon \left[ T_c - \Phi(x+k(x))(x+k(x)) \right].$$

We next show that $\mathbb{T}$ is a contraction mapping in the complete set $\mathcal{B}$ by estimating $\| \mathbb{T}k_1-\mathbb{T}k_2\|_{\mathcal{B}}$. Let  $k_1$ and $k_2$ in $\mathcal{B}$ whose norms are no more than $L$ be given. \\

\noindent To simplify the presentation we let $\xi=x+k_1(x)$ and $\zeta=x+k_2(x)$. 
\begin{align*}  
\|T k_1-T k_2 \|_{\mathcal{B}} &= \sup_{x \in \mathbb{R}} \epsilon | \left[ Tc-\Phi(x+k_1(x))(x+k_1(x))\right] -\left[ Tc-\Phi(x+k_2(x))(x+k_2(x))\right] |\\
& \le \epsilon \sup_{x \in \mathbb{R}} \left|  \Phi(\xi)(\xi) - \Phi(\zeta)(\zeta)  \right|\\
& \le \epsilon \left[ \sup_{x \in \mathbb{R}} |  \Phi(\xi)(\xi) - \Phi(\xi)(\zeta)|+ \sup_{x}|\Phi(\xi)(\zeta) - \Phi(\zeta)(\zeta)| \right] \\
& \le \epsilon \left[ L \|\xi-\zeta \|_{\mathcal{B}} +  \| \Phi(\xi)-\Phi(\zeta) \|_{\mathcal{B}} \right]  \text{, since } Lip_{\mathcal{B}}\Phi(\xi) \le L \\
& \le \epsilon  2L \| \xi - \zeta \|_{\mathcal{B}} \text{, since } \Phi \in \mathcal{G}_L \\
&\le \epsilon 2L \|k_1(x) - k_2(x) \|_{\mathcal{B}}.
\end{align*}

We thus showed that $\mathbb{T}$ is a contraction mapping. By the Banach contraction mapping theorem, there exists a unique $k^* \in \mathcal{B}$,  $k^*$ depends on $\Phi$, such that $\mathbb{T}k^*=k^*$. Letting $\xi = \eta+k^*(\eta)$ in \eqref{xi} finishes the proof. We call $\xi$ the \emph{preimage} of $\eta$ associated with $\Phi$.

\hfill \end{proof}

We now define the Budyko graph transform operator.\\
\begin{definition}\label{Bgto}
Let $\Phi \in \mathcal{G}_L$. We define the operator $\mathbb{M}$ acting on a graph $\Phi$ by: 
$$ \mathbb{M}(\Phi)(\eta)= \Phi(\xi)+F(\Phi(\xi), \xi)$$
where $\xi$  is the preimage of $\eta$ associated with $\Phi$ defined as in Lemma \ref{preimage} and $F$ as in the Definition \eqref{budxt}. We call $\mathbb{M}$ the Budyko graph transform operator.\\
\end{definition}

Before we show the contraction property of the Budyko graph transform operator, we proof the following computational Lemma that provides an upper bound estimate of $\epsilon$. \\

\begin{lemma} \label{epsize}
 Suppose that $F$ in  Definition \eqref{budxt} satisfies the Budyko parameter condition and the constant $L$ satisfies the estimate in Definition \ref{L}. Denote by $$ B_0 := 1-\frac{B}{R} \text{, } \delta:= \frac{L}{R} \text{ and } \kappa := \frac{B}{R(\delta+2B_0L)}.$$
If \begin{equation} \label{epsilon}\epsilon< \min \left\lbrace \frac{1}{8L},  \frac{\kappa}{1+4L\kappa} \right\rbrace  \end{equation} then
\begin{equation} \label{epsineq} 0<B_0 \left( 1+2 \rho L  \right)+\rho  \delta < 1 \end{equation}
with $\rho := \frac{\epsilon}{1-4 \epsilon L}$.\\
\end{lemma}

\begin{proof} 
Notice that since $\epsilon < \frac{1}{8L}$, by direct computation one may show that $0< \rho < 1$. Also, by the definition of the constant $L$ in Definition \ref{L}, $L > A \gg 2B$, one can show that $\delta+2B_0L = \frac{1}{R}\left( L-2B \right)+2L \gg 1$, and therefore, $0< \kappa < 1$. Since $0<\epsilon<  \frac{\kappa}{1+4L\kappa}$  then $\epsilon< \kappa \left( 1-4 \epsilon L \right).$ Also, notice that we may also write  $\kappa=\frac{1-B_0}{R(\delta+2B_0L)}$. Then the following sequence of inequalities leads to the claimed estimate \eqref{epsineq}:
\begin{align*}
0< \epsilon < \kappa (1-4 \epsilon L) & \Rightarrow 0< \frac{\epsilon}{1-4 \epsilon L} < \frac{1-B_0}{\delta+2B_0 L} \\
& \Rightarrow 0< \rho < \frac{1- B_0}{\delta + 2B_0 L} \\
& \Rightarrow 0<2B_0 L \rho + \rho \delta < 1-B_0 \\
& \Rightarrow 0< B_0(1+2\rho L)+\rho \delta  < 1.
\end{align*}
\hfill \end{proof}

\noindent We will use the same notation used in Lemma \eqref{epsize} for the remainder of this article. We now state the contraction result:\\
\begin{lemma}\label{main}
If the function $F$ in the extended Budyko map satisfies the Budyko parameters condition and $\epsilon$ satisfies  estimate \eqref{epsilon}, then the Budyko graph transform operator is a contraction mapping on the subset $\mathcal{G}_L$ of the space of graphs $\mathcal{G}$.\\
\end{lemma}

\begin{proof} 
First, we show that the operator $\mathbb{M}$ maps $\mathcal{G}_L$ to itself. Suppose that $\Phi \in \mathcal{G}_L$, and $\Phi$ can be written as $\Psi-T^*$ for some $\Psi \in \mathcal{G}$. Then for a given $\eta$, there exists $\xi$ such that:

\begin{align}\label{mphi}
\notag \mathbb{M}(\Phi)(\eta)  &= \Phi(\xi)+F(\Phi(\xi), \xi) \\
&= \left(\Psi - T^*\right)(\xi) - \frac{B+C}{R} \Psi(\xi)+\frac{C}{R}\overline{\Psi(\xi)}\\
\notag &\phantom{=} \left(\text{since } F(T^*(\xi),\xi)=0 \right)  \\
 &= \left( 1-\frac{B+C}{R}\right) (\Psi-T^*)(\xi)-\frac{B}{R}T^*(\xi)+\frac{C}{R} \left( \overline{\Psi(\xi)}-T^*(\xi) \right). 
\end{align}

\noindent Therefore, $\|\mathbb{M}(\Phi)\|_{\mathcal{G}} \le L$, since for any given $\eta$
\begin{align*}
\|\mathbb{M}(\Phi)(\eta)\|_{\mathcal{B}} & \le   \left(   1-\frac{1}{R}(B+C)\right)\|\Psi(\xi) - T^* (\xi)\|_\mathcal{B}+\frac{B}{R} \|T^*(\xi) \|_\mathcal{B} \\
& +\frac{C}{R} \|\Psi(\xi) - T^*(\xi) \|_\mathcal{B}  \\
& \le L.
\end{align*}

From equation \eqref{mphi} we may also estimate the Lipschitz constant of $M(\Phi)(\eta)$ for each fixed $\eta$. Recall that the Lipschitz constant of a function is unchanged upon adding a constant. Also, the Lipschitz constant of the sum of two Lipschitz functions is the sum of their Lipschitz constants. Therefore,

\begin{align*}
Lip_{\mathcal{B}}(\mathbb{M}(\Phi)(\eta))&=Lip_{\mathcal{B}}\left( \left(\Psi - T^*\right)(\xi) - \frac{B+C}{R} \Psi(\xi)+\frac{C}{R}\overline{\Psi(\xi)}  \right)\\
&\le \left( 1-\frac{B+C}{R}\right)Lip_{\mathcal{B}} (\Psi-T^*)(\xi) + Lip_{\mathcal{B}} \left(\frac{B+C}{R} T^*(\xi)\right) \\
&\le L.
\end{align*}

Next we show that $Lip_{\mathcal{G}}(\mathbb{M}(\Phi)) \le L$. Again, by rewriting $\mathbb{M}(\Phi)$ as in equation \eqref{mphi}, we have that 

\begin{align*}
Lip_{\mathcal{G}} \left( M(\Phi) \right) &=Lip_{\mathcal{G}} \left( \Phi -\frac{B+C}{R} \Psi  \right)\\
&\le  \left( 1-\frac{B+C}{R} \right)Lip_{\mathcal{G}} \left( \Phi \right) + \frac{B+C}{R}  Lip_{\mathcal{G}} \left( T^* \right)
&\le L.
\end{align*}
  
We have shown that $\mathbb{M}(\Phi)$  satisfies the norm condition and both the Lipschitz conditions of $\mathcal{G}_L$, and so we conclude $\mathbb{M}:\mathcal{G}_L \rightarrow \mathcal{G}_L$.  \\

Let $\Phi, \Psi \in \mathcal{G}_L$ and let $\eta \in \mathbb{R}$ be given. By the Preimage Lemma (\ref{preimage}), we have the preimages $x$ and $z$  associated with $\Phi$ and $\Psi$ such that
\begin{align*}
\eta &= x+\epsilon [\Phi(x)(x)-T_c] \\
\eta &= z+\epsilon [\Psi(z)(z)-T_c]. 
\end{align*}

We first compare $x$ and  $z$ to their respective graphs using the same notation as in Lemma (\ref{epsize}). Recall that $$\rho = \frac{\epsilon}{1-\epsilon 4L}. $$ Therefore, \\
\begin{align*}
|x-z| &= \epsilon |\Phi(x)(x) - \Psi(z)(z)| \\
 &\le \epsilon |\Phi(x)(x)-\Phi(x)(z)|+|\Phi(x)(z)-\Psi(x)(z)|+|\Psi(x)(z)-\Psi(z)(z)|\\
 & \le \epsilon \left( 2L |x-z | + \|\Phi-\Psi\|_{\mathcal{G}}  + 2L|x-z|+ \right)\\
 &\le \epsilon \left(  \|\Phi-\Psi \|_{\mathcal{G}} + \phantom{\int}4L |x-z| \right).
 \end{align*}
Therefore, \begin{equation}\label{xminz} |x-z| \le \rho \| \Phi - \Psi \|_{\mathcal{G}}. \end{equation} 

Another estimate yields 
\begin{align} \label{PhPs}
\|\Phi(x)-\Psi(z) \|_{\mathcal{B}} &\le \|\Phi(x)-\Psi(x) \|_{\mathcal{B}}+ \|\Psi(x) -\Psi(z) \|_{\mathcal{B}} \\ \notag
&\le \|\Phi-\Psi \|_{\mathcal{G}} + 2L|x-z |\\ \notag
&\le (1+\rho 2L) \|\Phi - \Psi \|_{\mathcal{G}}.
\end{align}

We now estimate the Budyko graph transform operator on the space of graphs $\mathcal{G}$. Fixing $\eta$, we have

$\| \mathbb{M}(\Phi)(\eta)-\mathbb{M}(\Psi)(\eta)\|_{\mathcal{B}} $
\begin{align}\label{Mcontr}
\notag &=\| \Phi(x)+F(\Phi(x),x) - \left(\Psi(z) + F(\Psi,z)\right)  \|_{\mathcal{B}}\\
\notag &= \| (1-\frac{1}{R}(B+C))\left[ \Phi(x)-\Psi(z) \right] + \frac{1}{R}Qs(y)\left[ \alpha(z)(y) -\alpha(x)(y) \right] \\
\notag &+\frac{1}{R} C \left[  \int_0^1 \Phi(x)(y) dy - \int_0^1 \Psi(z)(y) dy\right] \|_{\mathcal{B}}\\
\notag &\le \left(  |1-\frac{1}{R}(B+C)|+\frac{1}{R} C \right) \|\Phi(x)-\Psi(z) \|_{\mathcal{B}} + \frac{1}{R} L |x-z| \\
\notag &\le \left( 1-\frac{1}{R}B \right) \|\Phi(x)-\Psi(z) \|_{\mathcal{B}} + \delta \rho  \|\Phi -\Psi \|_{\mathcal{G}} \\
& \le \left( B_0(1+\rho 2L)+\rho \delta \right) \| \Phi -\Psi \|_{\mathcal{G}}.
\end{align}
 
By inequality \eqref{epsineq} in Lemma \ref{epsize}, we have showed here that the Budyko graph transform operator $\mathbb{M}$ is a contraction mapping on the complete set $\mathcal{G}_L$. Therefore, by the Banach contraction mapping theorem, we conclude that: \emph{There exists a unique Lipschitz continuous graph $\Phi^* \in \mathcal{G}_L$ such that $\mathbb{M}(\Phi^*)=\Phi^*$}. \\

The contraction constant in inequality \eqref{Mcontr} will reappear in the following proofs. It is then useful to denote this constant here. Let  \begin{equation}\label{gamma} \gamma = \left( B_0(1+\rho 2L)+\rho \delta \right).\\ \end{equation}

We finish the proof of Theorem \ref{mainresult} by showing the invariant and the attracting property of the Budyko map $m$ on the manifold  $\Phi^*(\mathbb{R}) \times \mathbb{R}$.
Recall from Lemma (\ref{preimage}), there exists a real valued function $k_{\Phi^*}$ such that $k_{\Phi^*}(x)=\epsilon\left( T_c- \Phi^*(x)(x)\right)$. Given $x \in \mathbb{R}$, let $$y=x-k_{\Phi^*}(x).$$ Then $$\Phi^*(y)=\mathbb{M}(\Phi^*(x)) \text{ and } y=x+\epsilon \left( \Phi^*(x)(x)-T_c \right).$$

\noindent And therefore, 
\begin{align}\label{mphistar}
m \left(\left[ \Phi^*(x),x \right] \right) &= \left[ \Phi^*(x)+F(\Phi^*(x),x), x+G(\Phi^*(x),x)\right] \notag \\ 
&= \left[ \mathbb{M}(\Phi^*)(x),x+G(\Phi^*(x),x)\right] \notag \\
&= \left[ \Phi^*(y), y\right].
\end{align}
\newline

So, $\Phi^*(\mathbb{R}) \times \mathbb{R}$ is invariant under the action of the Budyko map $m$, and thus, we have shown the invariant part, ie. part 1, of Theorem \ref{invariant}.  \\

We now show part 2 of Theorem \ref{invariant}, ie. show that the invariant manifold $\Phi^*(\mathbb{R}) \times \mathbb{R}$ is locally attracting. \\

Let $$\mathcal{C}:=\left\lbrace T \in \mathcal{B}: \|T \|_{\mathcal{B}} \le L \text{, } Lip_{\mathcal{B}} (T)\le L \right\rbrace.$$

We now show the contraction property of the Budyko map $m$ on $\mathcal{C} \times \mathbb{R}$ and thereby showing the attracting property of the invariant manifold. Given $T_0 \in \mathcal{C}$, then $T_0 = \Phi_0(\eta)$, where $\Phi_0(\eta) = T_0$ for any $\eta$. It is straightforward to check that $\Phi_0 \in \mathcal{G}_L$. Again, let $x \in \mathbb{R}$ be given, and let $y = x-k_{\Phi^*}(x)$ and $z= x-k_{\Phi_0}(x)$. Recall that both the constants $\gamma$  in equation \eqref{gamma} and $\rho$ in Lemma \ref{epsize} are positive and are majorized by 1.

\begin{align} \label{contraction}
\|  m([\Phi_0 , x])-m([\Phi^* , x]) \|_{\mathcal{B} \times \mathbb{R}} 
&=\|\left[\mathbb{ M }(\Phi_0)(z)-\mathbb{M}(\Phi^*)(y),z-y \right] \|_{\mathcal{B} \times \mathbb{R}} \notag \\
& \le \max \left\lbrace \gamma \|\Phi_0 - \Phi^*\|_{\mathcal{G}}, \rho\|\Phi_0 - \Phi^* \|_{\mathcal{G}} \right\rbrace  \\
& < \|[\Phi_0(x) - \Phi^*(x), x-x] \|_{\mathcal{B} \times \mathbb{R}}. \notag
\end{align}

Therefore, $\|m^k \left(T_0,x \right)-(\Phi^*,x) \| \rightarrow 0$ exponentially, thus, we have shown here part 2 of Theorem \ref{invariant}.
\hfill \end{proof}\\

Now that we have a result in  the embedded space $\mathcal{B} \times \mathbb{R}$, ie. there is a one dimensional locally attracting  manifold invariant under the extended Budyko map $m$, we will show the main result, which is the restriction of Theorem \ref{invariant} to the unit interval. \\

Recall that the goal of main result is to show the existence of the locally attracting invariant manifold of the Budyko map over the unit interval $I$. Suppose that $\eta, y \in I$, let $$\phi^*(\eta)(y) := \Phi^*(\eta)(y),$$ ie. the restriction of 
$\Phi^*(\eta)(y)$ over the unit interval. We will show that $$\mathcal{M}_I:=\left\lbrace [\phi^*(\eta)(y), \eta]:y, \eta \in I=[0,1] \right\rbrace $$ is the desired manifold.\\

We define a closed set  analogous to $\mathcal{C}$,$$\mathcal{C}_I= \left\lbrace T \in \mathcal{B}_I, \|T \|_{\mathcal{B}_I}\le L, Lip_{\mathcal{B}_I}(T) \le L \right\rbrace .$$

For the remainder of this proof, let $y, \eta \in I=[0,1]$. To show the first part of Theorem \ref{mainresult}, we show the contraction result \eqref{contraction} is also true in this restricted manifold under the restricted Budyko mapping $\widehat{m}$. Suppose that $\widehat{T}_0 \in \mathcal{C}_I$ and $\eta_0 \in I$, we extend $\widehat{T}_0$ in a similar fashion as the Budyko map $\widehat{m}$ and define $T_0(y) \in \mathcal{B}_I$
\begin{equation*}
T_0(y):=
\begin{cases}
&T_0(0) \text { ,  when } y<0\\
&T_0(y) \text{ ,  when } 0 \leq y \leq 1\\
&T_0(1) \text{ ,  when } 1<y.\\
\end{cases}
\end{equation*} 
In fact, $T_0(y) \in \mathcal{C}$. When the domain of the temperature profile is restricted to the unit interval, the Budyko mapping $m$ is precisely $\widehat{m}$. Furthermore, since for a given  $\eta \in [0,1]$ we have:
 \begin{equation*}
\| \widehat{m} \left( [\widehat{T}_0, \eta_0] \right) - [\phi^*(\eta)(y), \eta] \|_{\mathcal{B}_I \times I} 
\le \| {m} \left( [T_0, \eta_0] \right) - [\Phi^*(\eta)(y), \eta]\|_{\mathcal{B} \times \mathbb{R}}
\end{equation*}
then, the contraction result \eqref{contraction} is true also for the restricted space $\mathcal{B}_I \times I$. Hence, $\mathcal{M}_I$ is locally attracting under the action of $\widehat{m}$. \\

To show the second part of Theorem \ref{mainresult}, recall the equalities in \eqref{mphistar} shows that 
$$m \left([\Phi^*(\eta), \eta] \right)=[\Phi^*(\xi),(\xi)]. $$ 
\noindent Since we already specified that $\eta \in [0,1]$, then $m \left([\Phi^*(\eta), \eta] \right)=\widehat{m}\left([\Phi^*(\eta), \eta] \right)$, and if $\xi \in [0,1]$, then the first invariant alternative of the main result follows. Otherwise, $\xi \notin [0,1]$,  $\Phi^*(\xi)(y) \in \mathcal B$, so that the restriction of $\Phi^*(\xi)$ on the unit interval is also bounded continuous, hence in $\mathcal{B}_I$,  and the second alternative follows. This ends the proof of the main result Theorem \ref{mainresult}. \hfill \end{proof} \\
 
We now present the proof of Corollaries \ref{critset} and \ref{stability}.\\

\noindent \begin{proof} \textbf{(of Corollary \ref{critset})}\\

Again, we show the result holds in the embedded space, then show that the same result is also true in the restricted space. Let $\eta \in \mathbb{R}$. Since $\mathbb{M}(\Phi^*)(\mathbb R)=\Phi^*(\mathbb R)$, then the following holds:
\begin{align*}
\Phi^*(\eta) &= \Phi^*(\xi)+ F \left( \left[ \Phi^*(\xi)-T^*(\xi)+T^*(\xi), \xi \right] \right)
\end{align*}
\noindent where $\xi = \eta+k_{\Phi^{*}}(\eta)$.\\

Using reverse triangle inequality and the fact that $F(T^*(\xi),\xi)=0$ we find:
\begin{align*}
|| \Phi^*(\eta)-\Phi^*(\xi)||_{\mathcal{B}} &= || (B+C)\cdot [ \Phi^{*}(\xi) -T^*(\xi)] - C (\overline{\Phi^*(\xi)}-\overline{T^*(\xi)})||_{\mathcal{B}} \\
&\ge  \left|  (B+C)\cdot || \Phi^{*}(\xi) -T^*(\xi)||_{\mathcal{B}}- C || \overline{\Phi^*(\xi)}-\overline{T^*(\xi)}||_{\mathcal{B}}   \right| \\
& \ge B || \Phi^*(\xi) - T^*(\xi)||_{\mathcal{B}}.   
\end{align*}
\noindent In the last step of the above estimate, we used the fact that $$ \left|  \overline{\Phi^*(\xi)}-\overline{T^*(\xi)} \right| \le \| \Phi^*(\xi)-T^*(\xi)\|_{\mathcal{B}}.$$\\

\noindent We arrive at the desired estimate:

\begin{align} \label{oepsilon}
||\Phi^*(\xi)-T^*(\xi)||_{\mathcal{B}}& \le \frac{1}{B} ||\Phi^*(\eta)-\Phi^*(\eta + k_\Phi^*)\|_{\mathcal{B}} \notag \\
&\le \frac{1}{B} \cdot  \| k_{\Phi^*}(\eta)\| \notag \\
&\le \frac{\epsilon L}{B}.
\end{align} 

\noindent Since the estimate is true for any given ice line $\eta$ and since the upper bound $\frac{\epsilon L}{B}$ is independent of $\eta$, we conclude that $\|\Phi^*-T^*\|_{\mathcal{G}}< \frac{\epsilon L}{B}$. In particular, in the restricted space, ie. if $y, \xi \in [0,1]$, $\phi^*(\xi)(y)=\Phi^*(\xi)(y)$, and  the following inequality  holds $$ ||\phi^*(\xi)(y)-T^*(\xi)(y)||_{\mathcal{B}_I}\le ||\Phi^*(\xi)(y)-T^*(\xi)(y)||_{\mathcal{B}}.$$
Thus, the estimate in \eqref{oepsilon} applies to the restricted space, $\mathcal{B}_I$. The invariant manifold $\mathcal{M}_I$ is therefore within $O(\epsilon)$ of the critical set $\mathcal T_I^*:=\left\lbrace [T^*(\eta)(y),\eta]: \eta, y \in I=[0,1] \right\rbrace$.  This completes the proof of Corollary \ref{critset}.
 \hfill \end{proof}\\

\noindent \begin{proof} \textbf{of Corollary (\ref{stability})}\\
As an immediate result of Corollary \ref{critset}, we may estimate, within an order of $\epsilon$, the ice line temperature difference upperbound. That is: $$|\phi^*(\eta)(\eta)-T^*(\eta)(\eta)| < \frac{\epsilon L}{B}.$$  So that for a suficiently small $\epsilon$, $\phi^*(1)(1)<T_c$, prompting a condition for ice to form and ice line to move equatorward. Similar calculations yield the stability of the small ice cap and the snowball Earth and the instability of the big ice cap. The small ice cap is therefore a sink, while the big ice cap is a saddle. \hfill \end{proof}\\

\begin{figure}
\centering
\includegraphics[height=3in]{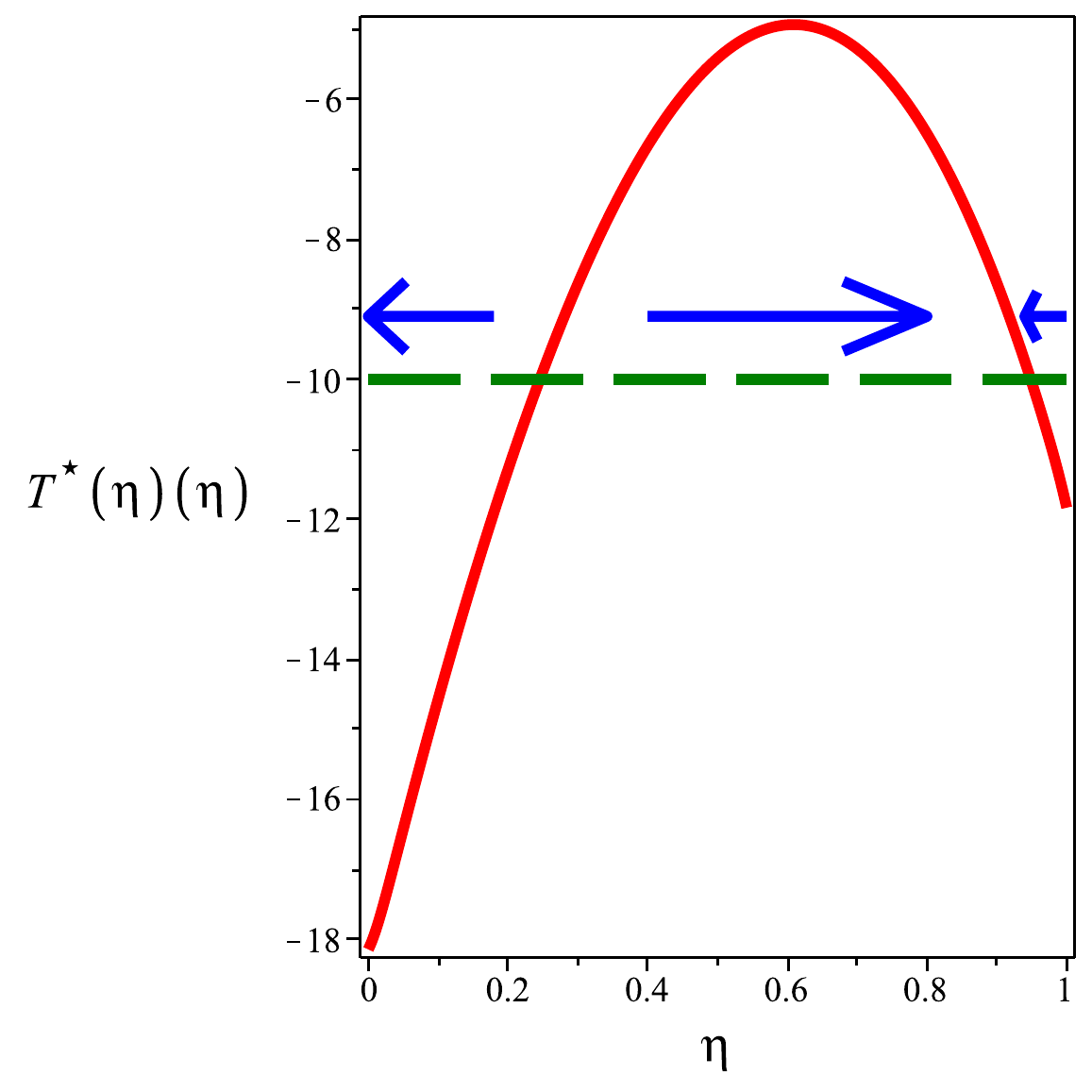}
\caption{The ice line temperature graph $T^*(\eta)(\eta)$ (in solid red) and the critical temperature $T_c=-10$ (in the green dash line) graphed using parameters  in \cite{tung07} over the  physically relevant interval  $[0,1]$. The arrows indicate the direction of the ice line movement. The intersections of the two curves $\eta \cong 0.2 \text{ and } 0.962$ are the large ice cap and the small ice cap equilibria. } \label{Tetaeta}
\end{figure}

Figure \ref{Tetaeta} illustrates the one dimensional dynamics of the slow system in Definition \eqref{budxt} dictated by  $G(T, \eta)=\epsilon \left( T(\eta)-T_c \right)$, using the parameters in Table \ref{parametertable}. If an ice line is initiated  equatorward of the large ice cap equilibrium, the ice line temperature of the equilibrium temperature profile  would be colder than the critical temperature $T_c$ and the ice line advances toward the equator leaving us with a snowball Earth. If an ice line is initiated between the large ice cap and the small ice cap, the ice line temperature would be warmer than the critical temperature and thus it would retreat  to the small ice cap equilibrium. Poleward of the small ice cap equilibrium, the ice line temperature is colder than the critical temperature causing ice to grow, hence, the ice line would advance toward the equator. It can be seen from this figure that the dynamics of the system would stay the same  for a value  $T_c$ in the range between about $-5.5^oC$ to about $-12^oC$, so long as the line $T_c$ cuts the solid (red) curve $T^*(\eta)(\eta)$ in two sites representing the large and the small ice cap equilibria.\\

\section{Bifurcation Diagrams and Animations}\label{anim}
\subsection{Bifurcation Diagrams}

\begin{figure}[h]
\centering
\includegraphics[width=4in]{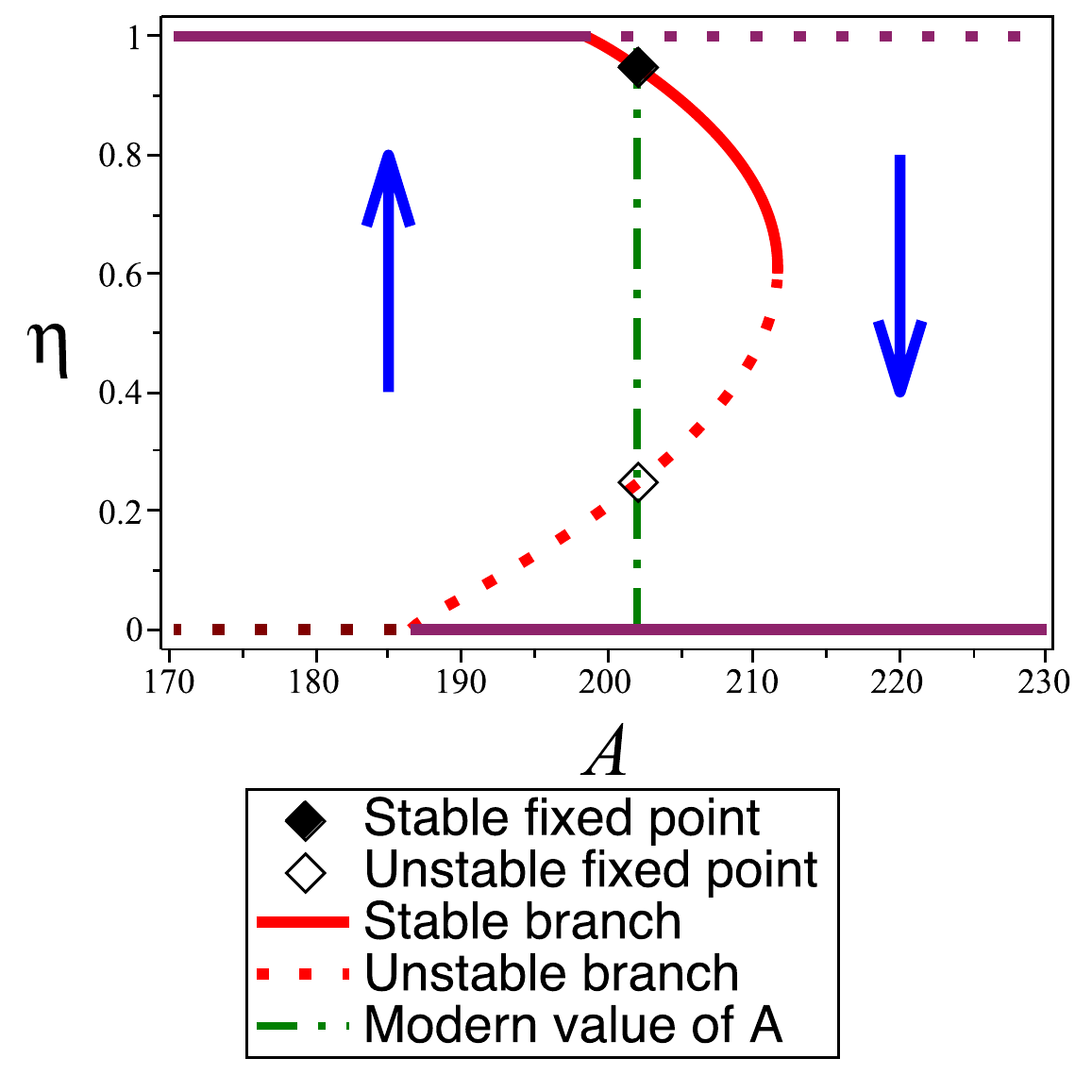}
\caption{The bifurcation diagram of the parameter $A$. All other parameters are as in Table \ref{parametertable}.}\label{Abif01}
\end{figure}

\begin{figure}
\centering
\includegraphics[width=4in]{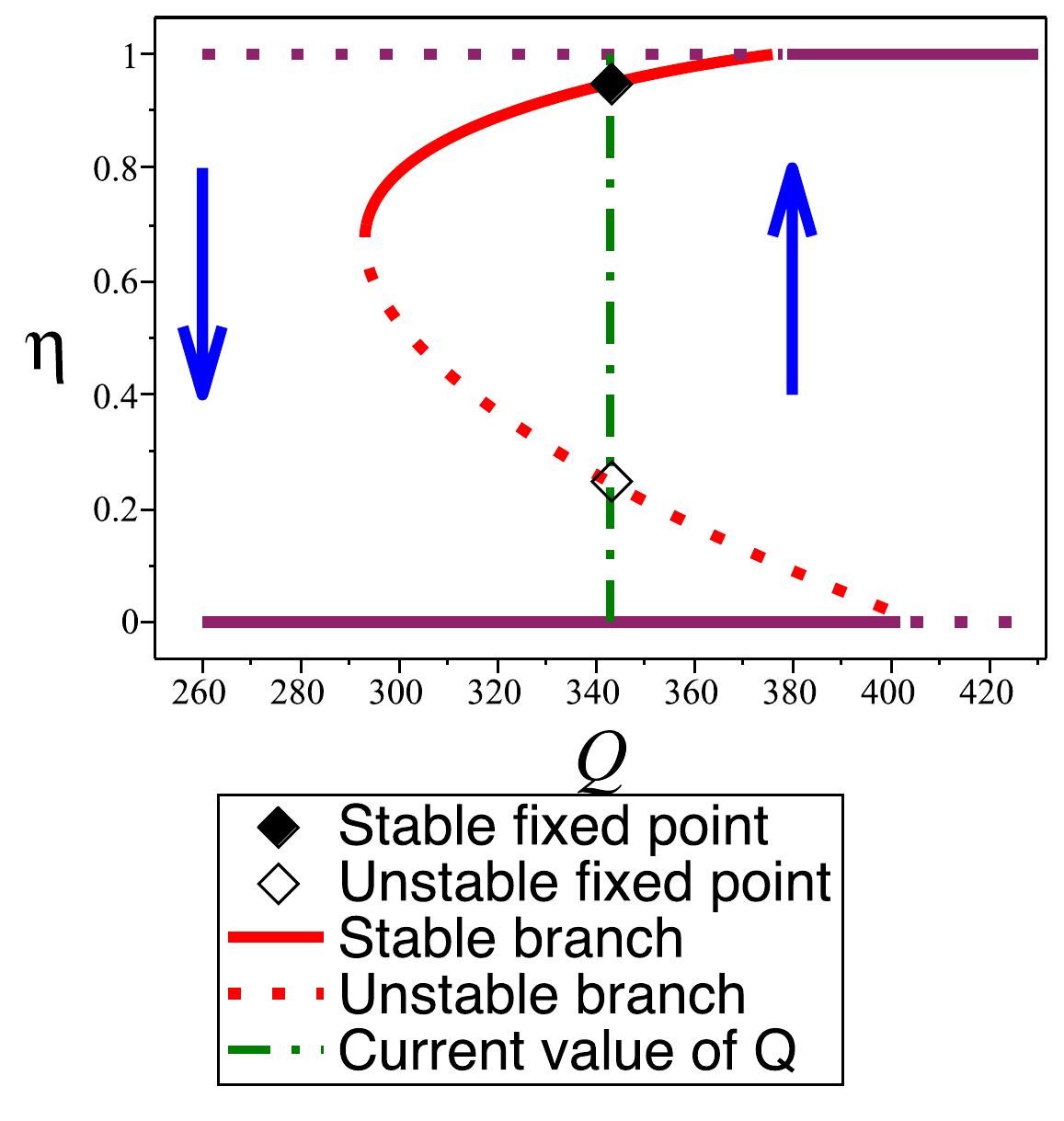}
\caption{The bifurcation diagram of the parameter $Q$. All other parameters are as in Table \ref{parametertable}.}\label{Qbif01}
\end{figure}


Readers may find bifurcation diagrams  similar  to Figures \ref{Qbif01} and \ref{Abif01} in other literatures (see \cite{dorian10, tung07}). We emphasize here that, unlike those presented in previous works, our stability analysis is a result of the ice line dynamics formulated in the slow equation, ie. \eqref{budt} of the fast-slow system in the statement of the main result Theorem \ref{mainresult}. This is indeed possible because the dynamics of the coupled system collapsed into a one dimensional system of the ice line.  \\

The bifurcation diagram of the system with respect to the \emph{greenhouse gas} parameter $A$ is presented in Figure \ref{Abif01}. The solid (red) curve in both figures is obtained by solving for the parameter $A$ when $T^*(\eta,\eta)$ is set to the critical temperature $T_c=-10^oC$ and the large green arrows indicate the dynamics of the ice line. If a pair of a parameter $A$ and an initial ice line $\eta$  is chosen to the right of this curve then the value of $T^*(\eta,\eta)$ for this pair will be less than the critical temperature $T_c$. In this case, the ice line equation dictates that the ice line moves down toward the equator, hence the downward blue arrow to the right of the red curve. Similarly, choosing a combination left of the red curve will move the ice line toward the pole.\\

The value $A=202$ in particular is chosen  in \cite{tung07} to indicate a modern value of radiative forcing obtained through satelite observations \cite{graves93}. With this value of $A$, the reader will observe the instability of the big ice cap equilibrium at $\eta \cong 0.2$ and  the small ice cap equilibrium $\eta \cong 0.962$. Furthermore, one recognizes the stability of the branches above the saddle node bifurcation point for the values near $A=212$ and $\eta=0.6$, represented by the solid red curve, and the instability of the lower branch, represented by the dash red curve. \\ 

We now focus on the physically relevant bifurcation points at the boundaries, ie. when $A \cong 186$ (corresponding to $\eta=0$) and  $A \cong 198$ (corresponding to $\eta=1$), where the determination of the bifurcation type is not so straightforward. The physically relevant condition that the ice line is not admissible outside of the unit interval imposes the mathematical restriction on the  dynamics of the ice line at the  boundaries, ie. $\eta=0$ and $\eta=1$, namely that 

\begin{align} \label{restriction}  &\text{at the equator or } \eta = 0 \, \quad \Delta \eta = \max \left\lbrace 0, \epsilon(T(\eta)(\eta)-T_c) \right\rbrace\\
\notag &\text{at the pole or } \eta = 1 \, \quad \Delta \eta = \min \left\lbrace 0, \epsilon(T(\eta)(\eta)-T_c) \right\rbrace. 
\end{align}

\noindent Using this restriction, one may determine the stability of branches in the bifurcation diagram at the physical boundary. In Figure \ref{Abif01}, the maroon solid lines at the  equator $\eta=0$  and the pole $\eta=1$ indicates a stable branch of the bifurcation diagram, suggesting an eternally ice covered or ice free state, depending on the parameter $A$. The maroon dash lines in the same figure indicates an unstable branch of the bifurcation diagram for $A$. \\

The reader may observe that the mathematical restrictions on $\Delta \eta$ in equation \eqref{restriction} introduces a non-differentiability of the ice line governing equation on a subset of the phase space. Such non-differentiability have known to introduce phenomena called \emph{border collision bifurcations} \cite{NuOttYorke}  or \emph{discontinuity induced bifurcation} \cite{koglen}. This suggests a more rigorous treatment of the Budyko-ice line model on the physical boundary, however, such treatment is beyond the scope of this paper.\\



A similar analysis is repeated for the bifurcation diagram of the parameter $Q$ in Figure \ref{Qbif01}. The phase line diagram for the current value of solar insolation $Q=343$ used in \cite{tung07} is shown in Figure \ref{Qbif01}, with blue arrows indicating the ice line dynamics. As in the discussion for the bifurcation diagram of parameter $A$, the stability of each branch of the bifurcation diagram for $Q$ is determined using the ice line equation \eqref{budt}, for ice line inside the unit interval and using \eqref{restriction} for the physical boundaries at the pole and the equator. The color scheme and the line style used in these figures are the same as that in Figure \ref{Abif01}.

\subsection{Animations}

\begin{figure}[h]
  \centering
   \subfloat{\label{strticov}} \href{run:staticicecovered.avi}{\includegraphics[width=1.5in]{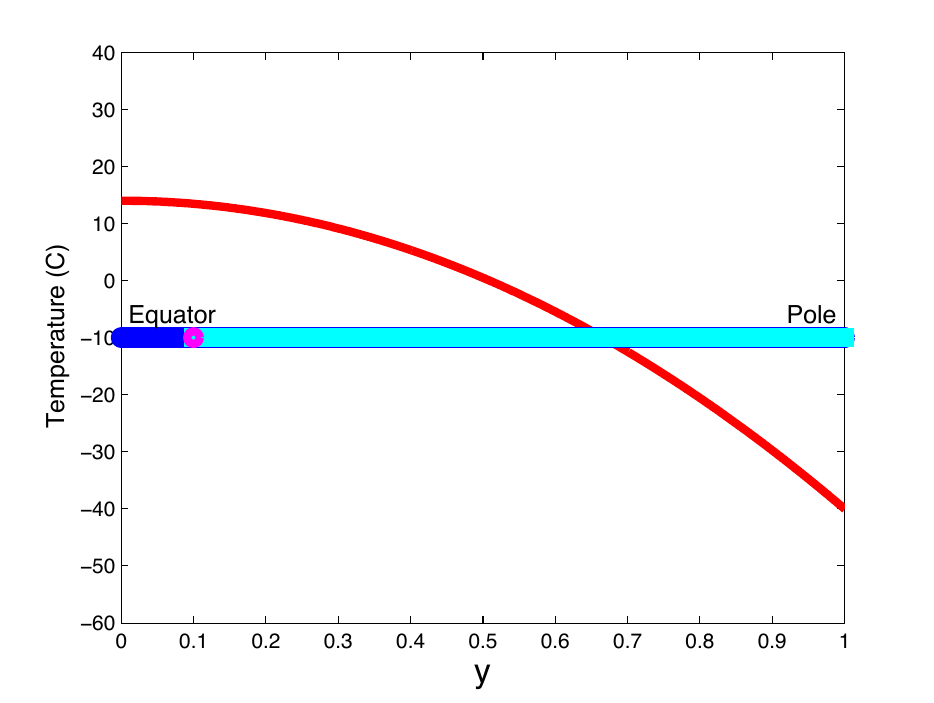}}                
   \subfloat{\label{strtmid}} \href{run:staticicemiddle.avi}{\includegraphics[width=1.5in]{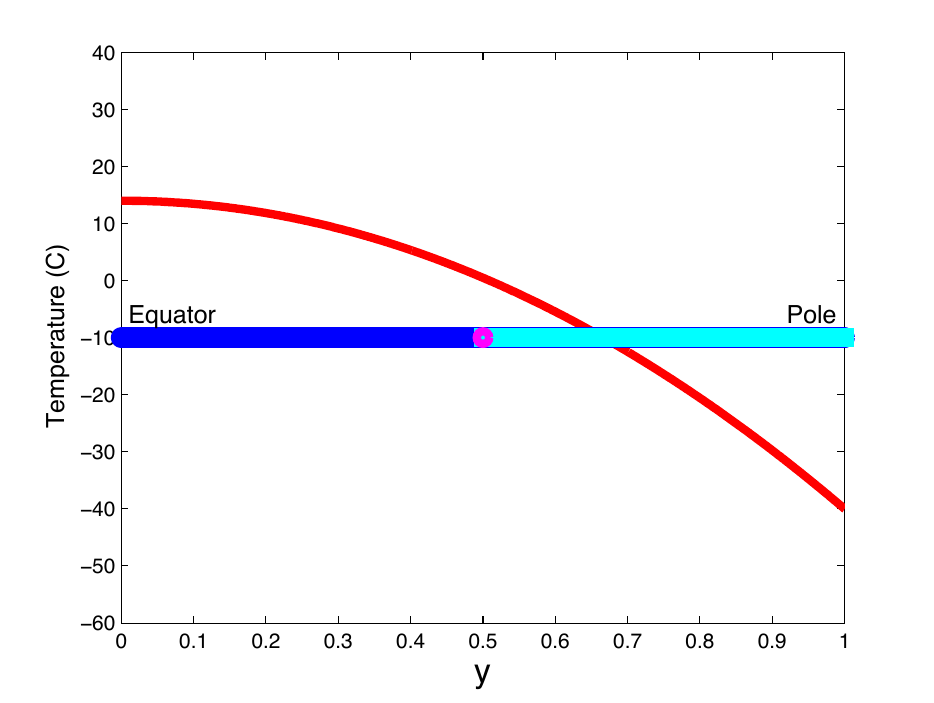}}   
   \subfloat{\label{strticf}} \href{run:staticicefree.avi}{\includegraphics[width=1.5in]{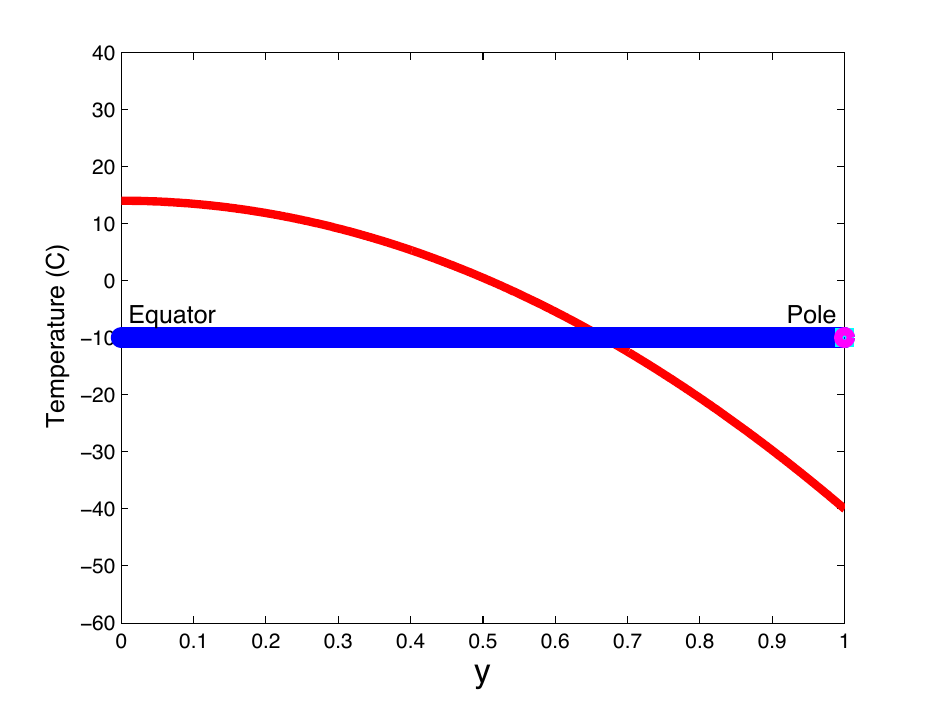}}  
   \subfloat{\label{endicov}\includegraphics[width=1.5in]{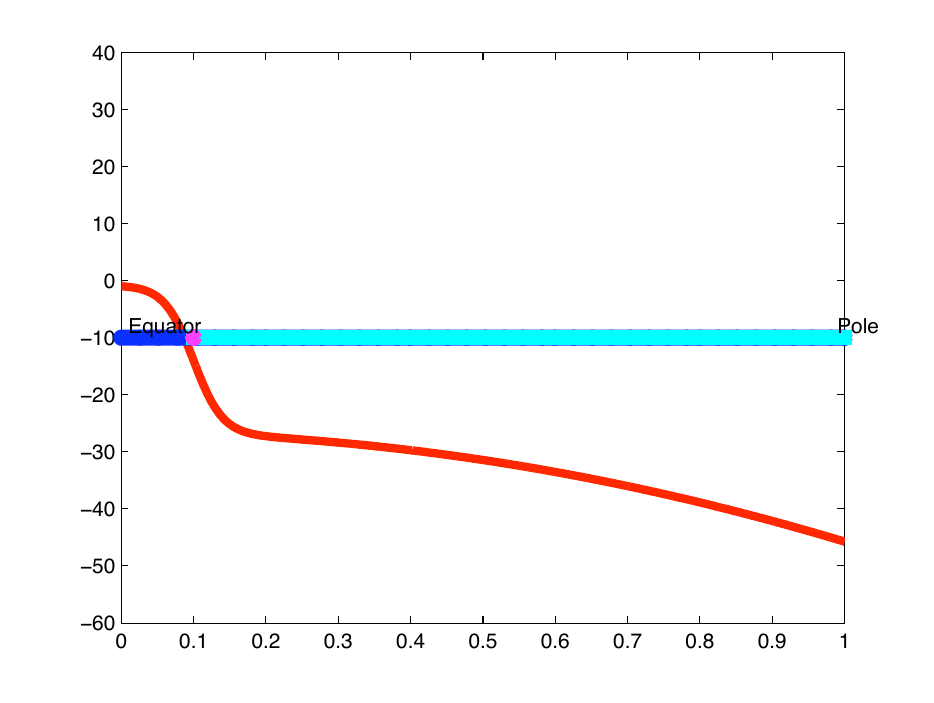}}
   \subfloat{\label{endmid}\includegraphics[width=1.5in]{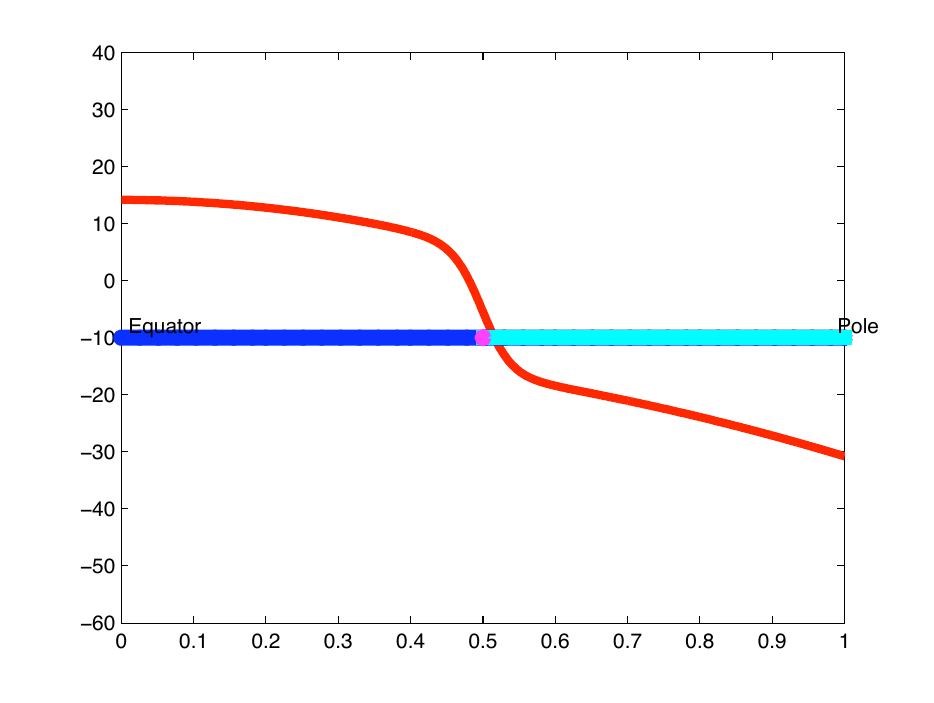}}  
   \subfloat{\label{endicf}\includegraphics[width=1.5in]{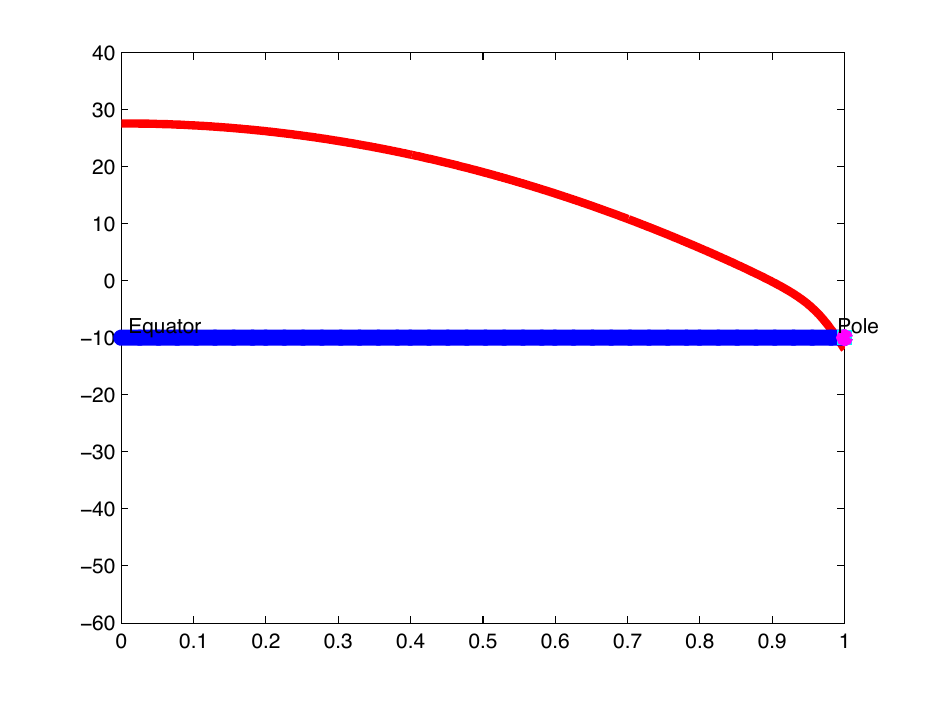}}
 
  \caption{\textbf{WITHOUT ICE LINE DYNAMICS}. The ice free surface is dark blue and the ice covered surface is bright cyan, while the red curve is the temperature profile. The vertical axes are the temperature in $Celcius$, the horizontal axes are the unit interval representing the sin of the latitude. The figures on the top are the initial temperature profile, $T(y)=14-54y^2$ with ice lines at $[0.1,T_c]$ (left), $[0.5,T_c]$ (middle) and $[1,T_c]$ (right), with $T_c=-10^oC$. The figures on the bottom are the steady state temperature profiles $T^*(\eta)(y)$ at the respective ice lines. Clicking on the top figures will start the animation. Adobe Acrobat Reader may be required.} 
\label{icf}
 \end{figure}

The animations (in Figure \ref{icf}) follow the evolution of the system without the ice line equation. In Figure \ref{icf} the red curve represents the temperature profile of the system, the light blue represents the ice covered region, while the dark blue represents the ice free region. The ice line of the boundary of these region is represented by a purple dot. The top row Figure \ref{icf} are the initial temperature profies, ie. the graph of $T(y)=14-54y^2$. The initial temperature profile  is chosen to represent a decreasing function with temperature at the equator at $14^oC$ and at the pole $-40^oC$ resembling a current average temperature at those extreme places. The parameters used in this animation can be found in Table \ref{parametertable}, which is based on parameter values used in \cite{tung07}. The bottom row of Figure \ref{icf} gives the equilibrium temperature profiles.

\begin{table} [h] 
\begin{center}
\addtolength{\tabcolsep}{1mm}
\renewcommand{\arraystretch}{1.2}
\begin{tabular}{|c|c|c|c|}

\hline
\textbf{Parameter }  & \textbf{Value} & \textbf{Unit} & \textbf{Source} \\
\hline
\hline 
Q & 343 & $Watts/m^2$  & \cite{tung07} \\
\hline 
$\alpha_1$ & 0.32 & Dimensionless  & \cite{tung07} \\
\hline 
$\alpha_2$ & 0.62 & Dimensionless   & \cite{tung07} \\
\hline 
A & 202 & $Watts/m^2$  & \cite{tung07} \\
\hline 
B & 1.9 & $Watts/(m^2 \phantom{0}^oC)$  & \cite{tung07} \\
\hline
C & 3.04 & $Watts/(m^2 \phantom{0}^oC)$  & \cite{tung07} \\
\hline
$T_c$ & -10 & $^oC$  & \cite{tung07}\\
\hline
$R$ & $12.6$ & $Watts/(m^2 \phantom{0}^oC)$  & \cite{mcgw}  \\
\hline
$M$ & 25 & none  & author\\
\hline
$\epsilon$ & $10^{-2}$ (dynamic) & $/^oC\text{year} $  & author \\
\hline 
$\epsilon$ & $0$ (static) & $/^oC\text{year} $  & author \\
\hline

\end{tabular}
\caption{The parameters used to create the simulations in Figures \ref{icf} and \ref{dyn}} \label{parametertable}

\end{center} 
\end{table}

\begin{figure}[h]
  \centering
  \subfloat{\label{strticecov}} \href{run:ebmsnowball.avi}{\includegraphics[width=1.5in]{static-icecovered.pdf}}                
  \subfloat{\label{strtic}} \href{run:ebmmiddle.avi}{\includegraphics[width=1.5in]{staticmiddle-start.pdf}}                
  \subfloat{\label{strticfree}} \href{run:ebmicefree.avi}{\includegraphics[width=1.5in]{staticicefree-start.pdf}}                 
  \subfloat{\label{endicecov}\includegraphics[width=1.5in]{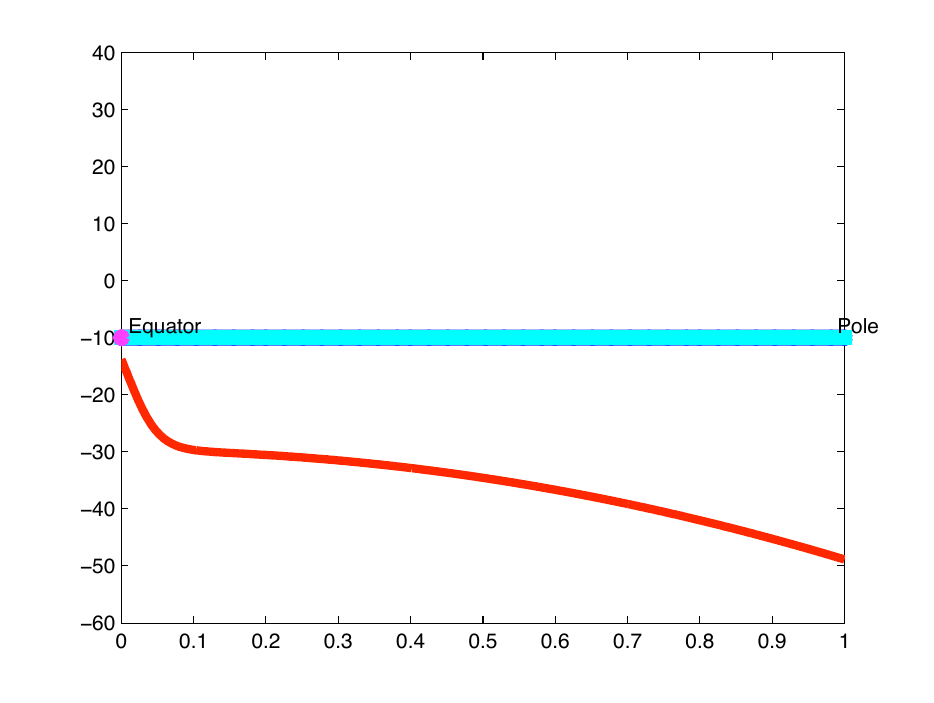}}  
  \subfloat{\label{endic }\includegraphics[width=1.5in]{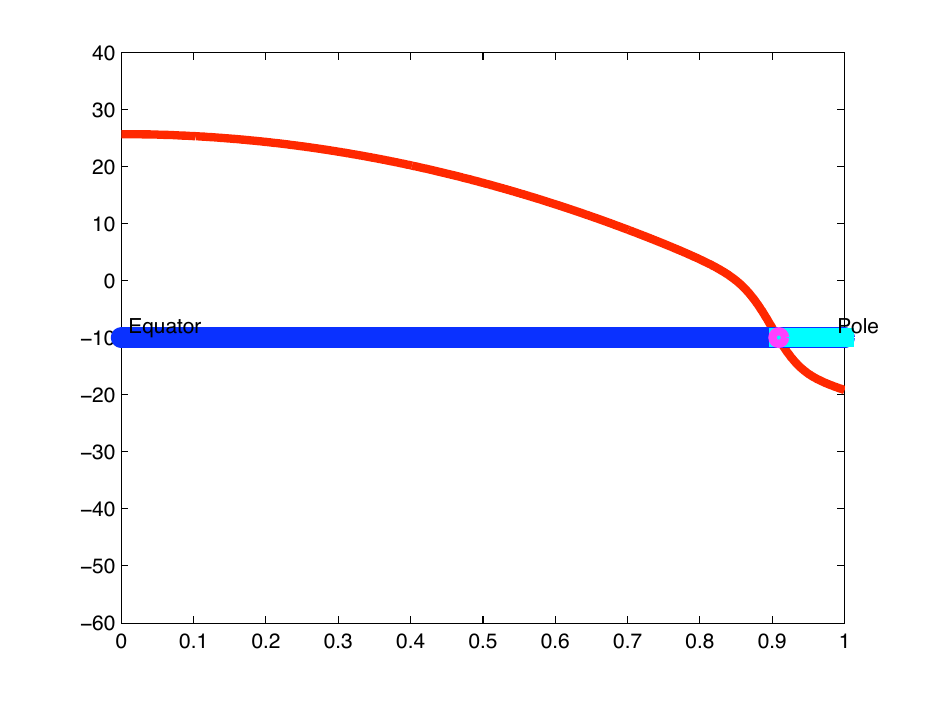}}
  \subfloat{\label{endicfree }\includegraphics[width=1.5in]{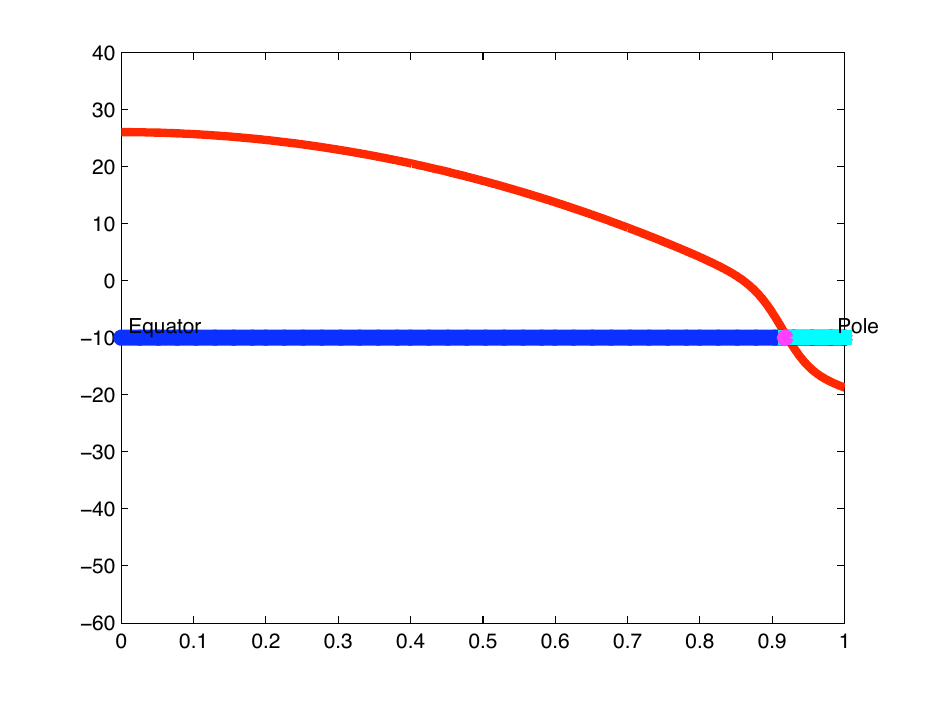}}
  \caption{\textbf{WITH ICE LINE DYNAMICS}. The figures on the top are the initial temperature profiles, $T(y)=14-54y^2$ with ice lines at $\eta_0=$ $0.1$, $0.5$ and $1$,  and the figure on the bottom are the steady state temperature profiles $T^*(\eta,y)$, with the ice lines at the equator or $\eta=0$, and $\eta \cong 0.962$. Again, the ice free surface is dark blue and the ice covered surface is bright cyan, while the red curve is the temperature profile. Clicking on the top figures will start the animations. Adobe Acrobat Reader may be required.}
   \label{dyn}
  \end{figure}


The  animations in Figure \ref{icf} illustrate the lack of ice line dynamics when the ice line equation is not included. In  Figure \ref{dyn}, the ice line equation is now included  and is coupled to the Budyko equation. Notice that in contrast to 
Figure \ref{icf}, the position of the ice line changes. The ice line dynamics can be observed in the animations in Figure \ref{dyn}.  We choose a value of $\epsilon = 10^{-2}$ that makes the animation run time reasonable. \\

\section{Concluding Discussion}\label{dir}
\subsection{Conclusion}
We introduce an auxilliary equation governing the ice line dynamics which is coupled to an ice albedo feedback EBM. The reformulation of the model contains two components: one equation representing planetary energy balance based on Budyko's EBM, and the other equation representing the dynamics of the ice line. This forms a system of fast-slow integro-difference equations and we show the existence of a one dimensional  attracting invariant manifold in this system. As a corollary of the main result, under parameter conditions prevalent in current literature, one may deduce the stability of various states of the planet, eg. ice free planet is unstable while a small ice cap is stable, etc. \\

\subsection{Future Directions} \label{futwk}
There are  several directions that one can explore based on this model. One immediate extension is to compute the invariant manifold explicitly similar to the work done by Foias, Sell and Titi \cite{fst}. The model formulation would be clarified even further with an explicit dynamics of  the ice line when it is at the equator or  the pole, possibly by introducing a discontinuity in the ice line equation. Another immediate improvement on this model is an extension to the southern hemisphere with two ice lines and a non symmetric transport. As stated in \cite{archer}, the two main feedbacks that dictate the interglaciation/ glaciation cycles are the ice albedo and greenhouse gasses feedbacks. The greenhouse gas effect in this model enters  through  the linear term $A+BT(y)$ affecting the radiative forcings. The work by Hogg in \cite{hogg} relates the evolution of temperature with that of carbondioxide as a response to the solar input variations caused by the Milankovitch cycle. A coupling of Budyko and Hogg's model may be fruitful in obtaining a simple model to understand the interglaciation/ glaciation cycles. The authors in \cite{n75}, \cite{n79}, \cite{n84} explore a similar model with only a diffusion transport. Another interesting modification would be to  include both diffusion (a second derivative term) and averaging (an integral term) in the transport process.\\

\section{Acknowledgement}
The author would like to thank Prof. R. McGehee, for his guidance in writing the PhD thesis that makes up the backbone of the main result and for his constant support during the course of writing this article. Also, the author is grateful for the support from the colleagues at the Mathematics Climate Research Network, as well as for many helpful suggestions from the editors and the anonymous referees. 
 
\bibliography{estherbib}
\bibliographystyle{siam} 
 
\end{document}